\documentclass[11pt,fleqn]{amsart} 

\usepackage{amsthm}
\usepackage{amsfonts}
\usepackage[english]{babel}
\usepackage{graphicx}
\usepackage{soul}
\usepackage{stfloats}
\usepackage{morefloats}
\usepackage{cite}
\usepackage{lscape}
\usepackage{epstopdf}
\usepackage{braket}
\usepackage[lite]{amsrefs}
\usepackage{mathbbol}
\usepackage{tikz,tikz-cd}
\usepackage{shuffle}

\usepackage{algorithm,algorithmicx,algpseudocode}

\usepackage{amsmath,amstext,amsopn,amsfonts,eucal,amssymb}
\usepackage{graphicx,wrapfig,url}

\newcommand\N{{\mathbb N}}

\newcommand\Q{{\mathbb Q}}
\newcommand\R{{\mathbb R}}

\newtheorem{theorem}{Theorem}[section]
\newtheorem{corollary}[theorem]{Corollary}
\newtheorem{lemma}[theorem]{Lemma}
\newtheorem{proposition}[theorem]{Proposition}
\newtheorem{definition}[theorem]{Definition}

\newtheorem{remark}[theorem]{Remark}

\newtheorem{criterion}[theorem]{Criterion}

\DeclareMathOperator{\Tr}{Tr}
\DeclareMathOperator{\Span}{span}
\DeclareMathOperator{\conv}{conv}
\DeclareMathOperator{\supp}{supp}

\begin{document}

	\title[Generalized Lefschetz Degree]{A Generalized Lefschetz Degree with applications to Fixed Point Theory}

	\author{Emanuele Zappala} 
	\address{Department of Mathematics and Statistics, Idaho State University\\
		Physical Science Complex |  921 S. 8th Ave., Stop 8085 | Pocatello, ID 83209} 
	\email{emanuelezappala@isu.edu}

	\maketitle

	\begin{abstract}
			We introduce and study a class of spaces, named Lefschetz admissible, where the identity maps on compacts can be approximated by certain maps with values in Lefschetz spaces in a suitable sense. 
			For maps defined on these spaces, we introduce a generalized Lefschetz degree obtained as the set of accumulation points of usual Lefschetz degrees for the approximant maps in a certain compactification. We show that when the degree is nonempty, fixed point results hold. We consider both the situation of spaces embedded in topological vector spaces (TVSs) as global and as neighborhood retracts. As two consequences, we obtain a generalization of a fixed point result due to Browder, and a no-retraction criterion for  compact subsets of locally convex TVSs. 
	\end{abstract}

	\date{\empty}

	\tableofcontents

	\section{Introduction}
	
	The study of fixed point theorems in infinite dimensional manifolds, in particular Banach spaces, has important applications in the study of integral equations and boundary value problems \cite{Amster}. Work of Browder \cite{Browder} has played a fundamental role in establishing early results in the field. In particular, Browder has obtained in \cite{Browder} some applications of the fixed point theorem of Lefschetz to the case of retracts in an infinite dimensional manifold. On the other hand, further work generalizing the Brouwer fixed point theorem by means of a generalized topological degree obtained via approximations, such as in \cite{Okon,Kaballo}, is applicable in topological vector spaces (TVSs) with the Klee admissiblity property. 
	
	The main objective of the present paper is to introduce a Lefschetz degree version of the definitions of \cite{Okon} for the Brouwer degree. This leads us to the notion of Lefschetz admissible TVSs, which we will show generalizes Klee admissibility. As an application, we will obtain a version of Browder's results to the case of the corresponding generalized Lefschetz degree. While Browder's theorems in \cite{Browder} require the space to be a Lefschetz space embedded in an infinite dimensional manifold as a (neighborhood) retract, in our setting the Lefschetz assumption falls only on the admissibility of the underlying TVS as an approximation property of sorts. In addition, while the TVSs in \cite{Browder} need to be quasi-complete and metrizable, the present setting allows us to discard these additional assumptions. 
	
	Similarly to the case of \cite{Okon}, which is the inspiring prototype for the present work, the notion of approximation compact maps is of central importance in guaranteeing the existence of a fixed point. In the Klee admissible case, work in \cite{Okon,Kim} has established classes of approximation compact maps for which fixed point theorems can be studied. Our main example of approximation compact map, in the generalized Lefschetz degree case, comes from compact maps, and it is of interest to obtain larger classes of such approximation compact maps. However, we leave this perspective to future work. 
	
	The present article is organized as follows. Section~\ref{sec:prelim} gives the fundamental definitions and notations used in the rest of the paper. Section~\ref{sec:TVS} introduces the case of generalized Lefschetz degree and Lefschetz admissibility, along with the associated fixed point properties, in the case of topological spaces that embed as rectracts in TVSs. In Section~\ref{sec:metric}, this framework is generalized to spaces that embed in TVSs as neighborhood retracts. Section~\ref{sec:applications} gives the two main applications of this work, namely a generalization of Theorems~3 and~4 in \cite{Browder} and a non-retraction criterion in locally convex TVSs. 
	
	\section{Preliminaries}\label{sec:prelim}
	
	Throughout the article, map will indicate a continuous map, although in some cases the term ``continuous'' is not explicitly mentioned. By compact map, we mean a map whose image has compact closure. 
	
	Following \cite{Rudin}, we recall the definition of topological vector space, which will be important throughout the paper. 
	\begin{definition}
		{\rm 
			A \emph{topological vector space} (TVS) over a topological field $\mathbb k$ is a vector space endowed with a topology that makes the scalar multiplication and addition operations continuous. 
		}
	\end{definition}
	We will consider just the case where $\mathbb k = \R$, and all TVSs are Hausdorff. In fact, we recall that the latter condition is not a very restrictive assumption since if $X$ is a non-Hausdorff TVS, the quotient space $X/\overline{\{0\}}$ obtained by quotienting out the closure of the singleton containing the zero element of $X$ is a Hausdorff TVS. Over the reals (or complex numbers), a subset $A$ of a TVS is said to be \emph{balanced}, or \emph{circled},  if $\alpha A \subset A$ for every $|\alpha| \leq 1$. Moreover, it is known that every TVS has a fundamental system of open balanced neighborhoods of $0$, where a fundametnal system of $0$ is a set which is a local basis at $0$. For a given TVS, we will denote by $\mathcal U$ such system of open balanced neighborhoods, which we consider fixed once a TVS is given. In particular, this will allow  us to use that fact that if $x\in U\in \mathcal U$, then $-x\in U$ as well. This fact will be used multiple times with no further mention. 
	Morevoer, we recall that a TVS is said to be locally convex when such $\mathcal U$ can also be chosen to consist of convex sets. When the underlying TVS is locally convex, we will therefore assume that the chosen $\mathcal U$ consists of open balanced and convex sets. 
	
	Let $X$ denote a topological space. The \emph{Lebesgue covering dimension}, or \emph{topological dimension}, of $X$ is defined as the minimal natural number $n\in \N$, if it exists, such that every finite open covering of $X$ has an open refinement with multiplicity at most equal to $n+1$. When a finite $n$ cannot be found, the topological dimension is said to be infinite. We also recall that the open refinement can be assumed to be finite without changing the definition, since the terms in any open refinement can be grouped in such a way that one obtains a finite open refinement without changing the original multiplicity. 
	We point out some notational discrepancies in the literature that we need to settle in order to avoid any potential confusion in the following. We mention that some authors, e.g. \cite{Munkres}, allow arbitrary open coverings in the definition of Lebesgue covering dimension. The definition used in the present work is therefore less restrictive. Moreover, in traditional texts it is common to find the definition of Lebesgue covering dimension given only for normal spaces. This includes the influential reference \cite{Engelking}. As explained in \cite{Engelking}, the reason for such restriction does not really stem from the fact that for more general topological spaces the definition would not be meaningful, but rather from the fact that it otherwise ``presents some undesirable features of a separation axiom''. Some authors consider the most general case, which we follow in the present article, as for example seen in \cite{Katetov3}.  
	
	We recall the following result from dimension theory, found in \cite[Theorem 3.1.4]{Engelking}, where the notion of covering dimension for a normal space refers to finite coverings. 
	
	\begin{theorem}
			For every closed subspace $M$ of a normal space $X$ we have $\dim M \leq \dim X$. 
	\end{theorem}
	 With the definitions used in the present article, the same result remains valid also when $X$ is not assumed to be normal. However, we remark that other frequently found definitions of finite dimensionality fail to satisfy this property, and therefore this is very specific to the finite open cover version of the definition used here. 
	 Interestingly, by Dowker's theorem \cite[Theorem 3.2.1]{Engelking}, for a paracompact Hausdorff space the definition of finite covering dimension can be equivalently reformulated using arbitrary open covering. Throughout the article, we will emphasize the normality of the spaces when this implies no further restriction, e.g. when dealing with paracompact Hausdorff spaces. The reader who wants to follow Engelking's definition closely, might add normality as an underlying assumption when we discard it. This is remarked explicitly in Section~\ref{sec:metric}. 

	We recall that a TVS $X$ is said to be \emph{Klee admissible} when for each compact $K\subset X$ and fundamental neighborhood of zero $U\in \mathcal U$, we can find a continuous map $P_K : K \longrightarrow X$ such that $x-P_K(x)\in U$ for each $x\in K$ and $P_K$ is contained in a vector subspace of $X$ of finite rank. We remark that, as pointed out also in \cite{Okon}, this definition can be equivalently given with maps $P_K : X \longrightarrow X$, as a procedure akin to Tietze extension theorem allows to extend $P_K$ from $K$ to the whole $X$. When we mention Klee admissibility, we will assume that the definition with $P_K$ defined on the whole space $X$ is being considered.  
	
	Two of the most commonly used notions of subnet are due to Kelley \cite{Kelley}, and Willard \cite{Willard}. 
	In the present article we will consider the notion of subnet in Willard's sense. Let $\Gamma$ and $\Gamma'$ be directed sets. Given a  net $x: \Gamma \longrightarrow X$, we say that a net $y: \Gamma' \longrightarrow X$ is a subnet of $x$ if we can find a \emph{cofinal morphism} $\phi: \Gamma' \longrightarrow \Gamma$ such that $y = x\circ \phi$. Morphism here means that $\phi$ is order preserving, and being cofinal means that for every $\gamma \in \Gamma$ there exists $\gamma'_0 \in \Gamma'$ such that $\phi(\gamma'_0) \geq_\Gamma \gamma$. This condition on $\phi$ means that the image of $\phi$ is \emph{cofinal}. In particular, by the order preservation of $\phi$, we have that $\phi(\gamma') \geq_\Gamma \gamma$ for every $\gamma'\in \Gamma'$ with $\gamma' \geq_{\Gamma'} \gamma'_0$. For a net $\{x_\gamma\}_{\gamma\in\Gamma}$, we abbreviate the reindexed net $\{x_{\phi(\gamma')}\}_{\gamma'\in\Gamma'}$ by $\{x_{\gamma'}\}_{\gamma'\in\Gamma'}$. Thus $x_{\gamma'}$ stands for $x_{\phi(\gamma')}$, and convergence is taken with respect to the directed order on the domain $\Gamma'$ of $\phi$. To indicate the existence of a cofinal morphism $\phi$ and the corresponding subnet $\{x_{\phi(\gamma')}\}_{\gamma'\in\Gamma'}$, we also say that $\Gamma'$ is a subnet of $\Gamma$, which follows the same nomenclature adopted in \cite{Okon}. In what follows we omit the subscripts on the order symbols whenever the index set is clear from context.

	\section{Lefschetz admissible topological vector spaces and a topological degree}\label{sec:TVS}
		
	In this article we will call Lefschetz space, following \cite{Browder}, a topological space $X$ such that for any continuous self-mapping $T : X \longrightarrow X$, the Lefschetz degree $\Lambda(T)$ is well defined, and whenever $\Lambda(T) \neq 0$, $T$ has a fixed point. We recall that the Lefschetz degree $\Lambda(T)$ is defined as the alternating series
	\begin{eqnarray}\label{eqn:Lefschetz_deg}
			\Lambda(T) &=& \sum_{s=0}^\infty (-1)^s\Tr(T_s), 
	\end{eqnarray} 	
	where $T_s$ is the map induced by $T$ on singular homology $H_s(X,\Q)$, with coefficients in the rationals. Using the same convention of \cite{Browder}, $\Lambda(T)$ is said to be well defined when $T_s$ is zero except for possibly finitely many $s\in \N$, and  the trace of $T_s$ is well defined for all $s$. 
	
	Lefschetz spaces are homologically small, in the following sense. It was proved by Browder \cite[Lemma~2]{Browder} that if $Y$ is a Lefschetz space, then $H_s(Y,\Q)$ is finite dimensional for all $s$, and it is nontrivial for finitely many $s\in \N$. We recall the fact that compact-convex subsets of $\R^n$ are Lefschetz spaces. 
	
	The following lemma will be useful in the rest of the article, and it will be used repeatedly. 
	
	\begin{lemma}\label{lem:homeo_Lefschetz}
			If $X$ is homeomorphic to a Lefschetz space $Y$, then $X$ is a Lefschetz space as well. 
	\end{lemma}
	\begin{proof}
			Let $f : X \longrightarrow Y$ denote a homeomorphism. Let $T : X \longrightarrow X$ be a continuous self-mapping. Define the continuous map $\hat T : Y \longrightarrow Y$ by $\hat T = f\circ T\circ f^{-1}$. Since $f_* : H_*(X,\Q) \longrightarrow H_*(Y,\Q)$ is an isomorphism in homology, and since $\hat T_s$ is nontrivial only for finitely many $s\in \N$, $T_s$ is also nontrivial for finitely many $s\in \N$. Moreover, from
			\begin{eqnarray*}
					\Tr(\hat T_*) 
					&=& \Tr(f_*\circ T_*\circ f^{-1}_*)\\
					&=& \Tr(T_*), 
			\end{eqnarray*}
			it follows that the trace of $T_s$ is well defined for all $s\in \N$, and $\Lambda(T)$ is well defined. It also follows that $\Lambda(\hat T) = \Lambda(T)$, and if $\Lambda(T) \neq 0$, it follows that $\hat T$ has a fixed point $y\in Y$. Setting $x = f^{-1}(y)$ we see that $T(x) = x$, which concludes the proof. 
	\end{proof}

	In what follows we consider to work on a Hausdorff space $X$ which is a topological vector space (TVS) over $\R$ with a system of circled neighborhoods of zero denoted by $\mathcal U$.  We consider a notion of approximation net similar to the one considered in \cite{Okon}. 
	
	\begin{definition}\label{def:approx_net}
		{\rm 
				Let $X$ be a (Hausdorff) TVS, and let $\Gamma$ indicate a directed set. A family of continuous maps $\{P_\gamma : X \longrightarrow X_\gamma\}_{\gamma \in \Gamma}$ is said to be an {\it L-approximation net}, or simply approximation net, if $X_\gamma \subset X$ is a Lefschetz space for each $\gamma \in \Gamma$, and if for every compact set $K \subset X$ and neighbor of zero $U\in \mathcal U$, we can find $\gamma_0\in \Gamma$ such that
				\begin{eqnarray*}
						x-P_\gamma(x) &\in& U, 
				\end{eqnarray*}
				for all $x\in K$, whenever $\gamma \geq \gamma_0$. 
		}
	\end{definition} 

	The following will be useful.   
	\begin{lemma}\label{lem:approx_subnet} 
			Let $\{P_\gamma\}_{\gamma\in \Gamma}$ denote an L-approximation net, and let $\{P_{\gamma'}\}_{\gamma'\in \Gamma'}$ be a subnet. Then, $\{P_{\gamma'}\}_{\gamma'\in \Gamma'}$ is also an L-approximation net. 
	\end{lemma}

	\begin{definition}
		{\rm 
				We say that a Hausdorff TVS $X$ with a circled system of neighbors of zero $\mathcal U$ is a Lefschetz-admissible space, or L-admissible for short, if $X$ admits an L-approximation net.  
		}
	\end{definition}

	While \cite{Okon,Kim} consider Klee admissible spaces, with approximation nets having finite dimensional image, we are considering a similar situation where our target spaces are homologically finite dimensional, as observed above. 

	The following result constructs our prototype for Lefschetz admissible spaces. 
	
	\begin{proposition}\label{pro:normed_Ladmissible}
			Let $X$ be a normed space where $\mathcal U$ is the system of balls centered at zero. Then $X$ is L-admissible. 
	\end{proposition}
	\begin{proof}
			We introduce the following directed set: $\Gamma := \{(K,\epsilon)\ |\ \emptyset \neq K \subset X\ {\rm is\ compact},\ \epsilon>0 \}$, ordered by $(K_1,\epsilon_1) \geq (K_2,\epsilon_2)$ if and only if $K_1 \supset K_2$ and $\epsilon_1 \leq \epsilon_2$. 
			Given $\gamma := (K,\epsilon)\in \Gamma$, we want to define $P_\gamma : X \longrightarrow X_\gamma$, where $X_\gamma$ is a Lefschetz space. We find finitely many $x_i\in K$, with $i = 1, \ldots, n$ such that $K \subset \bigcup_{i=1}^nB(x_i,\epsilon)$. Let $\mathcal A := \{B(x_i,\epsilon)\}_{i=1}^n \cup \{X-K\}$ be a covering of $X$, and let $\{\eta_i\}_{i=1}^n\cup\{\eta_\infty\}$ be a partition of unity subordinate to $\mathcal A$, where each $\eta_i$ is supported inside $B(x_i,\epsilon)$ for $i=1,\ldots, n$, and $\eta_\infty$ is supported in $X\setminus K$. We choose a point $y \in X$ arbitrarily, and let
			$X_\gamma := \overline{\conv}\langle x_1, \ldots, x_n , y\rangle$ be the closed convex hull generated by $x_1, \ldots, x_n , y$. Since $X_\gamma$ is convex, it is contractible, and therefore any self-mapping of $X_\gamma$ has a well defined Lefschetz degree equal to $1$.  Moreover, $X_\gamma$ is linearly homeomorphic to a closed and bounded subset of $\R^m$ for some $m$, and it is therefore compact. By the Schauder fixed point theorem, any self-map of $X_\gamma$ has a fixed point. It therefore follows that $X_\gamma$ is a Lefschetz space. 
			We define $P_\gamma : X \longrightarrow X_\gamma$ as
			\begin{eqnarray*}
					P_\gamma(x) &:=& \sum_{i=1}^n \eta_i(x)x_i + \eta_\infty(x)y. 
			\end{eqnarray*} 
			Then, for any $x\in K$ we have
			\begin{eqnarray*}
					\|x-P_\gamma(x)\| &=& \|x - \sum_{i=1}^n \eta_i(x)x_i - \eta_\infty(x)y\|\\
					&=& \|\sum_{i=1}^n \eta_i(x)x -  \sum_{i=1}^n \eta_i(x)x_i\|\\
					&\leq& \sum_{i=1}^n \eta_i(x)\|x-x_i\|\\
					&<& \epsilon. 
			\end{eqnarray*}
			For each $\gamma \in \Gamma$ we construct such a $P_\gamma$, and we want to show that $\{P_\gamma\}_{\gamma\in \Gamma}$ is an L-approximation net. 
			Let $K\neq \emptyset$ be compact, and let $B(0,\epsilon)\in \mathcal U$. Set $\gamma_0 := (K,\epsilon)$, and let $\gamma = (K_1,\epsilon_1)$ be such that $\gamma \geq \gamma_0$. Then, by definition we have $K_1 \supset K$ and $\epsilon_1 \leq \epsilon$. Therefore, for each $x\in K \subset K_1$ we have
			\begin{eqnarray*}
					\|x-P_\gamma(x)\| &<& \epsilon_1 \leq \epsilon, 
			\end{eqnarray*}
			which shows that $x-P_\gamma(x)\in U$, therefore completing the proof. 
	\end{proof}

	By adjusting the proof of the previous result, we can show a more general version of it. 
	\begin{proposition}\label{pro:TVS_Ladmissible}
				Locally convex, paracompact Hausdorff TVSs are L-admissible spaces. 
	\end{proposition}
	\begin{proof}
			We define the directed set $\Gamma$ as the set of pairs $(K,U)$ where $K$ is a nonempty compact in $X$, and $U\in \mathcal U$, ordered as
			\begin{eqnarray*}
					(K_1,U_1) \geq (K_2,U_2) &\overset{\text{def}}{\iff}& K_1 \supset K_2 \ \land \ U_1 \subset U_2. 
			\end{eqnarray*}
			For a compact $K\neq \emptyset$ and a given $U\in \mathcal U$, we can cover $K$ by opens of type $x+U$, $x\in K$, and extract a finite subcovering $\{x_1+U, \ldots, x_n+U\}$. Choosing an arbitrary $y\in X$, we define the covering $\mathcal A = 	\{x_1+U, \ldots, x_n+U, X \setminus K\}$, and set $X_\gamma = \overline{\conv}\langle x_1, \ldots, x_n,y\rangle$ be the corresponding closed convex hull of $x_1, \ldots, x_n, y$, where $\gamma = (K,U) \in \Gamma$.  We let $\{\eta_1, \ldots, \eta_n, \eta_\infty\}$ denote a partition of unity subordinate to $\mathcal A$, where $\eta_i$ has support inside $x_i+ U$, and $\eta_\infty$ has support in $X\setminus K$. Then, we define $P_\gamma : X \longrightarrow X_\gamma$ by the assignment
			\begin{eqnarray*}
					P_\gamma(x) &=&  \sum_{i=1}^n \eta_i(x)x_i + \eta_\infty(x)y, 
			\end{eqnarray*}		
			which is continuous by construction. Let us denote by $V = \Span \{x_1, \ldots, x_n,y\}$, which is linearly homeomorphic to $\mathbb R^m$, where $m\in \mathbb N$ is the size of a maximal subset of $\{x_1,\ldots, x_n, y\}$ consisting of linearly independent vectors. We denote this map as $\phi$. Moreover, as seen for example in Theorem 1.21 in \cite{Rudin}, the linear homeomorphism is constructed by taking a basis of $V$, say $\{x_{i_1}, \ldots, x_{i_m}\}$, and mapping it to the standard basis of $\R^m$. It follows that $X_\gamma$ is mapped inside a closed ball of $\R^m$. Since homeomorphisms are closed maps, it follows that $X_\gamma$ is identified with a closed subset of a compact subset of $\R^m$, and therefore $\phi(X_\gamma)$ is compact in $\R^m$. From the fact that $X_\gamma$ is convex and $\phi$ is linear, it also follows that $\phi(X_\gamma)$ is convex. Since compact convex subsets of $\R^m$ are Lefschetz spaces, it follows that $X_\gamma$ is homeomorphic to a Lefschetz space, and by Lemma~\ref{lem:homeo_Lefschetz} $X_\gamma$ is a Lefschetz space too. To complete the proof, we need to show that we can approximate the identity over compacts of $X$. For a compact $K$ and a $U\in \mathcal U$, we set $\gamma_0 = (K,U)$. Then, for an arbitrary $x\in K$, we let $x_{j_1}, \ldots, x_{j_k}$ denote the elements of $x_1, \ldots, x_n$ such that $x\in x_{j_r} + U$. In particular, $\eta_i(x) =  0$ for $i \neq j_1, \ldots, j_k$, and clearly $\eta_\infty(x) = 0$. Then, we have
			\begin{equation*}
					x-P_\gamma(x) = x - \eta_{j_1}(x)x_{j_1} - \cdots - \eta_{j_k}(x)x_{j_k} = \eta_{j_1}(x)(x-x_{j_1}) + \cdots + \eta_{j_k}(x)(x-x_{j_k}) \in U, 
			\end{equation*}
			since $U$ is convex. 
			Then, for any $\gamma = (K_1,U_1)$ with $\gamma \geq \gamma_0$, we have that
			\begin{eqnarray*}
					x-P_\gamma(x) \in U_1, 
			\end{eqnarray*}
			for each $x\in K_1$, and by definition of ordering in $\Gamma$ it also follows that 
						\begin{eqnarray*}
				x-P_\gamma(x) \in U, 
			\end{eqnarray*}
			for each $x\in K$.
	\end{proof}

	More generally, we see that Klee-admissible TVSs are also L-admissible. 

	\begin{proposition}\label{pro:Klee_Ladmissible}
			Let $X$ be a Klee-admissible TVS. Then, $X$ is L-admissible. 
	\end{proposition}
	\begin{proof}
			As in the proof of Proposition~\ref{pro:TVS_Ladmissible}, we define the directed set $\Gamma$ of pairs $(K,U)$ where $K$ is nonempty and compact, and $U$ is a balanced neighborhood of zero, ordered by inclusion with respect to $K$, and reverse inclusion with respect to $U$. Given $\gamma := (K,U) \in \Gamma$, let $P_\gamma: X \longrightarrow X_\gamma$ denote the approximation of the identity on $K$ whose existence is guaranteed by Klee-admissibility, where $X_\gamma$ is a finite dimensional subspace of $X$. Let $\phi_\gamma: X_\gamma \longrightarrow \R^{n_\gamma}$ denote a linear homeomorphism of $X_\gamma$ with a Euclidean space. Consider the compact convex space $Y_\gamma$ obtained by taking the closed convex hull of the compact set $\phi_\gamma(P_\gamma(K))$. Since $Y_\gamma$ is contractible, the Lefschetz degree is always well defined, and moreover, by the Brouwer fixed point theorem a continuous self-map of $Y_\gamma$ has a fixed point. Therefore, $Y_\gamma$ is a Lefschetz space. Let $P'_\gamma := \phi_\gamma^{-1}r_\gamma\phi_\gamma P_\gamma : X \longrightarrow X'_\gamma$, where $r_\gamma$ is a retraction of $\R^{n_\gamma}$ onto the compact-convex subset $Y_\gamma$, and we have set $X'_\gamma := \phi_\gamma^{-1}(Y_\gamma)$. By Lemma~\ref{lem:homeo_Lefschetz}, we have that $X'_\gamma$ is a Lefschetz space. Applying this procedure for each $\gamma = (K,U)$ as above, we therefore construct a net $\{P'_\gamma\}_{\gamma\in \Gamma}$, which we now show to be an L-approximation net. Let $K$ be compact and $U\in \mathcal U$, set $\gamma = (K,U)$. By Klee-admissibility we have $x-P_\gamma(x)\in U$ for each $x\in K$. Then, we have
			\begin{equation*}
					x-P'_\gamma(x) = x - \phi_\gamma^{-1}r_\gamma\phi_\gamma P_\gamma(x) = x - \phi_\gamma^{-1}\phi_\gamma P_\gamma(x) = x - P_\gamma(x) \in U. 
			\end{equation*}
			As in the proof of Proposition~\ref{pro:TVS_Ladmissible} we see that for any $(K_1,U_1) = \gamma' \geq \gamma$ we have $x-P'_{\gamma'}(x) \in U_1$ for all $x\in K_1$, hence $x-P'_{\gamma'}(x) \in U$ for all $x\in K$, which completes the proof. 
	\end{proof}

	The following is a strengthening of the known fact that locally convex TVSs are Klee admissible. We provide a direct proof where the Klee approximant maps are uniformly continuous on the whole TVS via an extension result of Kat\v{e}tov. The proof will be useful for later results in Section~\ref{sec:metric}. 
	\begin{proposition}\label{pro:no_paracompact}
			Locally convex TVSs are Klee-admissible, where the approximant maps are uniformly continuous over the whole space. In particular, they are L-admissible. 
	\end{proposition}
	\begin{proof}
			Let $X$ be a locally convex TVS, let $K$ be a nonempty compact subset, and let $U\in \mathcal U$ be a convex fundamental neighborhood of zero. Since $X$ is a Hausdorff TVS, it is completely regular \cite{Kelley}, and we can therefore find for each $x\in K$ a function $u_x : X \longrightarrow [0,1]$ such that $u_x(x) = 1$ and $u_x(X\setminus \{x+U\}) = \{0\}$. We consider the open covering $\{u_x^{-1}((0,1])\ |\ x\in K \}$ and we let $x_1, \ldots, x_n \in K$ be such that, $\{u_i^{-1}((0,1])\ |\ i=1, \ldots, n \}$ is a subcovering of $K$, where we indicated $u_i$ the function corresponding to each $x_i$. We define the map $\eta_i : K \longrightarrow [0,1]$ by
			\begin{eqnarray*}
					\eta_i(x) &=& \frac{u_i(x)}{u_1(x)+ \cdots + u_n(x)}. 
			\end{eqnarray*}
			Applying Kat\v etov's results \cite{Katetov,Katetov2}, we can extend each $\eta_i$ to a uniformly continuous map $\hat \eta_i : X \longrightarrow [0,1]$. We define the map $P: X\longrightarrow X$ as 
			\begin{eqnarray*}
					P(x) &=& \sum_{i=1}^n \hat \eta_i(x)x_i.
 			\end{eqnarray*}
 			For $x\in K$, we have that $x\in u_{i_t}^{-1}((0,1])$ for $t = 1, \ldots, r$, and since $u_{i_t}(x) \neq 0$ implies that $x\in x_{i_t} + U$ by construction of $u_{i_t}$, it follows that $x-x_{i_t} \in U$ for each $t=1,\ldots, r$, while $u_i(x) = 0$ for each $i\neq i_t$ for all $t$. Then, we have
 			\begin{eqnarray*}
 					x - P(x) &=& \sum_{t=1}^r \eta_{i_t}(x)(x-x_{i_t})\in U. 
 			\end{eqnarray*} 
 			Since $P(X) \subset\Span \{x_1, \ldots, x_n \}$, $P$ has finite rank image. The second statement follows from Proposition~\ref{pro:Klee_Ladmissible}. 
	\end{proof}
	\begin{remark}
		{\rm
				In fact, the same construction of the proof of Proposition~\ref{pro:no_paracompact}, applying Kat\v{e}tov's extension theorem gives that a Klee admissible Hausdorff TVS admits approximants that are uniformly continuous and defined on the whole TVS. In other words, Klee admissibility can be stated equivalently as the existence, for each compact $K$ and $U\in \mathcal U$, of a uniformly continuous map $P: X \longrightarrow E$ with $E$ of finite rank, such that $x-P(x)\in U$ for each $x\in K$. This is a strengthening of the discussion given in Section~\ref{sec:prelim}. We have not been able to locate an explicit reference of this fact in previous literature. 
		}
	\end{remark}
	We point out that examples of Klee-admissible spaces without the local convexity assumption have been the object of active study, and explicit examples can be found for instance in \cite{Klee}, or also by combining results of \cite{Hadzic,Ishii}. Therefore, Proposition~\ref{pro:Klee_Ladmissible} gives examples of L-admissible spaces in the non-locally convex setting. 

	We now consider a notion of approximation compact map similar to \cite{Okon}, adapted to our setting. Since the version of \cite{Okon} is a generalization of the notion of A-proper map of \cite{Brow-Petr}, it follows that the approximation compact maps considered below are closely related to A-proper maps as well. As discussed in Section~\ref{sec:prelim}, we will simplify notation without explicitly reporting the cofinal morphisms $\phi: \Gamma' \longrightarrow \Gamma$ that define subnets. Therefore, $x_{\gamma'}$ will indicate $x_{\phi(\gamma')}$, $x_{\gamma''}$ will indicate $x_{\phi(\phi'(\gamma''))}$ and so on, leaving understood that there is an underlying $\phi$ or $\phi'$. We will say that $\Gamma'$ is a subnet of $\Gamma$ to indicate that there is a cofinal morphism function $\phi: \Gamma' \longrightarrow \Gamma$, where both $\Gamma'$ and $\Gamma$ are directed sets. 
	
	\begin{definition}\label{def:approx_comp_0}
			{\rm 
					Let $X$ be a Lefschetz admissible space, and let $\{P_\gamma\}_{\gamma\in\Gamma}$ denote an L-approximation net. Let $\emptyset \neq C \subset X$ be a closed subset of $X$ such that $X$ retracts onto $C$ via the retraction $r: X \longrightarrow C$. For each $\gamma \in \Gamma$ we let $r_{\gamma}$ denote $r_{|X_\gamma} : X_\gamma \longrightarrow C$. Let $T: C \longrightarrow X$ be a continuous map. Then, we say that $T$ is approximation compact with respect to $\{P_\gamma\}_{\gamma\in \Gamma}$ and $r$, or just approximation compact if no ambiguity arises, if the following condition is satisfied. For any subnet $\Gamma'$ of $\Gamma$, and sysetm $\{x_{\gamma'}\}_{\gamma'\in \Gamma'}$ such that $x_{\gamma'} \in X_{\gamma'}$ and $\lim_{\gamma'} (x_{\gamma'}-P_{\gamma'}Tr_{\gamma'}x_{\gamma'}) = 0$, there exist a subnet $\{x_{\gamma''}\}_{\gamma''\in \Gamma''}$ of $\{x_{\gamma'}\}_{\gamma'\in \Gamma'}$ and $x\in C$ such that $\lim_{\gamma''} x_{\gamma''} = x$ and $T(x) = x$. 
			}
	\end{definition}
	We can extend the previous definition to the case of fixed point equations of type $T(x) + y = x$ as follows. 
	
	\begin{definition}\label{def:approx_comp_Y}
			{\rm 
					Let $X$ and $T$ be as in Defintion~\ref{def:approx_comp_0}. For a subset $Y \subset X$, we say that $T$ is approximation compact with respect to $\{P_\gamma\}_{\gamma\in \Gamma}$, $Y$ and $r$, or just approximation compact if no ambiguity arises, if the map $T_y(x) := T(x) + y$ is approximation compact (in the sense of Definition~\ref{def:approx_comp_0}) with respect to  $\{P_\gamma\}_{\gamma\in \Gamma}$ and $r$ for every $y\in Y$. 
			}
	\end{definition}

	We note that setting $Y = \{0\}$ recovers Definition~\ref{def:approx_comp_0} from Definition~\ref{def:approx_comp_Y}. 
		\begin{lemma}\label{lem:y_shift_approx}
		Let $T : C \longrightarrow X$ be approximation compact with respect to $\{P_\gamma\}_{\gamma \in \Gamma}$, $Y \subset X$, and $r : X \longrightarrow C$. For a given $y\in Y$, define the map $T_y(x) := T(x) + y$. Then, $T_y$ is approximation compact with respect to $\{P_\gamma\}_{\gamma \in \Gamma}$, $\{0\}$, and $r$.
	\end{lemma}
 
	Our main example of approximation compact map now follows. 
	
	\begin{theorem}\label{thm:approx_compact_compact_im}
			Let $T : C \longrightarrow X$ be a compact continuous map, $X$ is L-admissible with an approximation net $\{P_\gamma\}_{\gamma\in \Gamma}$, and $r : X \longrightarrow C$ is a retraction. Let $Y$ denote the set of points of $X$ such that $T(C) + y \subset C$. Then, $T$ is approximation compact with respect to $\{P_\gamma\}_{\gamma\in \Gamma}$, $Y$ and $r$. 
	\end{theorem}
	\begin{proof}
				We consider the case $Y = \{0\}$ for simplicity, since we can otherwise just rename the map $T_y$ by $T$ and proceed analogously since $T_y$ satisfies the conditions in Definition~\ref{def:approx_comp_0}.   
				Let $\Gamma'$ be a subnet of $\Gamma$, and let $\{x_{\gamma'}\}_{\gamma'\in \Gamma'}$ be such that 
				\begin{eqnarray}\label{eqn:copmact_approx_1}
						\lim_{\gamma'} (x_{\gamma'}-P_{\gamma'}Tr_{\gamma'}x_{\gamma'}) = 0, 
				\end{eqnarray} 
				where $x_{\gamma'} \in X_{\gamma'}$ for each $\gamma' \in \Gamma'$. We want to show that we can extract a subnet of $\{x_{\gamma'}\}_{\gamma'\in \Gamma'}$ which converges to some $x\in C$ with $x-T(x) = 0$. 
				
				Let $W$ be a neighborhood of zero, and let $U\in \mathcal U$ be another neighborhood such that $U+U \subset W$, whose existence is guaranteed by the continuity of $+$ in a TVS. By \eqref{eqn:copmact_approx_1}, we can find $\gamma_0\in \Gamma'$ such that $\gamma' \in \Gamma'$ with $\gamma'\geq \gamma_0$ implies $x_{\gamma'} - P_{\gamma'}Tr_{\gamma'}x_{\gamma'} \in U$. Since $\overline{T(C)}$ is compact, and $\{P_{\gamma'}\}_{\gamma'\in\Gamma'}$ is an approximation net by Lemma~\ref{lem:approx_subnet}, we can find an index $\gamma_1\in \Gamma'$ such that $\Gamma'\ni\gamma'\geq \gamma_1$ implies 
				\begin{eqnarray*}
				Tr_{\gamma'} x_{\gamma'} - P_{\gamma'}Tr_{\gamma'} x_{\gamma'} \in U. 
				\end{eqnarray*} 
				Therefore, for $\Gamma' \ni \gamma_2 \geq \gamma_0, \gamma_1$, it follows that $\gamma'\in \Gamma'$ and $\gamma' \geq \gamma_2$ implies 
				\begin{eqnarray*}
						x_{\gamma'} - Tr_{\gamma'}x_{\gamma'} 
						&=&  x_{\gamma'} - P_{\gamma'}Tr_{\gamma'} x_{\gamma'} + P_{\gamma'}Tr_{\gamma'} x_{\gamma'} - Tr_{\gamma'}x_{\gamma'}
						 \in U + U \subset  W,  
				\end{eqnarray*}
				which shows that 
				\begin{eqnarray}\label{eqn:x_conv}
						\lim_{\gamma'} x_{\gamma'} - Tr_{\gamma'}x_{\gamma'} &=& 0. 
				\end{eqnarray}
				Since $\overline{T(C)}$ is compact, we can extract a subnet $x_{\gamma''}$, such that $Tr_{\gamma''}x_{\gamma''}$ converges.  
				Since $T(C) \subset C$, $C$ is closed (being a retract), and applying \eqref{eqn:x_conv} it follows that $\exists x:= \lim_{\gamma''} x_{\gamma''} \in C$. Then, we have
				\begin{equation*}
						T(x) = T(r(x)) = T(r(\lim_{\gamma''}x_{\gamma''})) = \lim_{\gamma''} T(r_{\gamma''}x_{\gamma''}) = \lim_{\gamma''} x_{\gamma''} = x. 
				\end{equation*}
	\end{proof}

				A concrete situation where Theorem~\ref{thm:approx_compact_compact_im} can be applied to obtain approximation compact operators is the following. Let $X$ be a normed space, along with an approximation net $\{P_\gamma\}_{\gamma\in\Gamma}$ as constructed for example in Proposition~\ref{pro:normed_Ladmissible}. Let $T : X \longrightarrow X$ be a completely continuous operator such that $T(\overline{B(0,R)}) \subset \overline{B(0,R)}$. Then, with the radial retraction $r: X \longrightarrow \overline{B(0,R)}$ and $Y = \{0\}$ we have that the restriction of $T$ to the closed ball $\overline{B(0,R)}$ respects the hypotheses of Theorem~\ref{thm:approx_compact_compact_im}, and it is therefore approximation compact. This situation arises frequently in the study of boundary value problems \cite{Amster}. When $\overline{T(\overline{B(0,R)})} \subset B(0,R)$, for $y\in B(0,r)$ for a sufficiently small $r$ we have the condition that $T(\overline{B(0,R)}) + y \subset \overline{B(0,R)}$. The same construction can be adjusted when the balls are not centered at zero. 

	In what follows, we use the extended real line $\R_\infty = \R\cup \{\pm \infty\}$, endowed with the topology generated by the basis $\mathcal B = \{(a,b)\ |\ a,b\in \R \} \cup \{[-\infty,a)\ |\ a\in \R \} \cup \{(a,+\infty]\ |\ a\in \R \}$. We notice that $\R_\infty$ is compact. 

	We introduce now a generalized version of the Lefschetz degree, in the same spirit of the construction of \cite{Okon} for the Brouwer degree. 
	
	\begin{definition}\label{def:gen_deg}
		{\rm 
				Let $T : C \longrightarrow X$ be a continuous map, where $X$ is L-admissible with L-approximation net $\{P_\gamma\}_{\gamma \in \Gamma}$, and let $r : X \longrightarrow C$ be a retraction. Then, we define the generalized Lefschetz degree of $T$, denoted by $\Lambda(T)$, as 
				\begin{eqnarray*}
						\Lambda(T) &=& \{\Lambda(P_\gamma Tr_\gamma)\ |\ \gamma\in \Gamma \}', 
				\end{eqnarray*}
				where $\Lambda(P_\gamma Tr_\gamma)$ is the Lefschetz degree of the map $P_\gamma Tr_\gamma : X_\gamma \longrightarrow X_\gamma$, and the prime symbol indicates the accumulation set of the indexed net in the topological space $\R_\infty$. In other words, $\Lambda(T)$ consists of the elements of $\R_\infty$ that are limits of subnets of $\{\Lambda(P_\gamma Tr_\gamma)\}_{\gamma\in \Gamma}$. 
		}
	\end{definition}

	\begin{lemma}\label{lem:deg_homotopy}
			Let $H : C \times [0,1] \longrightarrow X$ be a homotopy between $F$ and $G$, which are as in Definition~\ref{def:gen_deg}. Then $\Lambda(P_\gamma Fr_\gamma) = \Lambda(P_\gamma Gr_\gamma)$ for each $\gamma\in \Gamma$. In particular, $\Lambda(T) = \Lambda(G)$. 
	\end{lemma}
	\begin{proof}
			Let $H : C \times [0,1] \longrightarrow X$ be a homotopy between $F$ and $G$. Then, for each $\gamma\in \Gamma$, we have an induced homotopy $H_\gamma: X_\gamma \times [0,1] \longrightarrow X_\gamma$ defined as $H_\gamma(x,t) = P_\gamma H(r_\gamma(x),t)$, between $P_\gamma Fr_\gamma$ and $P_\gamma Gr_\gamma$, which shows that $\Lambda(P_\gamma Fr_\gamma) = \Lambda(P_\gamma Gr_\gamma)$. Therefore, we also have $\Lambda(F) = \Lambda(G)$. 
	\end{proof}

	\begin{theorem}\label{thm:fixed_point}
				Let $T : C \longrightarrow X$ be approximation compact with respect to $\{P_\gamma\}_{\gamma \in \Gamma}$, $Y \subset X$, and $r : X \longrightarrow C$. Then, the following facts hold. 
				\begin{itemize}
					\item[(i)] 
							We have $\Lambda(T) \neq \emptyset$. 
					\item[(ii)]
							Suppose $a\in \Lambda(T)$, and let $\{a_{\gamma'}\}_{\gamma'\in \Gamma'}$ be a subnet of $\{\Lambda(P_\gamma Tr_\gamma)\ |\ \gamma\in \Gamma \}$ such that  $\lim_{\gamma'} a_{\gamma'} = a$. Then, if $a\neq 0$ or $a = 0$ and $a_{\gamma'} \neq 0$ for all $\gamma' \geq \gamma'_0$ for some $\gamma'_0\in \Gamma'$, the fixed point equation $x - T(x) = y$ has a solution for every $y\in Y$.  
				\end{itemize}
	\end{theorem}
	\begin{proof}
				(i) We first show that $\{\Lambda(P_\gamma Tr_\gamma)\ |\ \gamma\in \Gamma \} \neq \emptyset$. This follows from the fact that $X_\gamma$ is a Lefschetz space for each $\gamma \in \Gamma$. In fact, the map $P_\gamma T r_\gamma : X_\gamma \longrightarrow X_\gamma$ is continuous, and therefore by definition of Lefschetz space we have that $\Lambda(P_\gamma Tr_\gamma) \in \Q$ is well defined. 
				Since $\R_\infty$ is compact, we can extract a converging subnet from $\{\Lambda(P_\gamma Tr_\gamma)\ |\ \gamma\in \Gamma \}$, which therefore has accumulation points in $\R_\infty$. 
				
				\noindent (ii) 
				 Let $a\in \Lambda(T)$, with $a_{\gamma'} \rightarrow a$, and $a_{\gamma'} = \Lambda(P_{\gamma'} Tr_{\gamma'})$ such that either $a\neq 0$, or $a_{\gamma'} \neq 0$ for each $\gamma' \geq \gamma'_0$. For $y\in Y$, the map $T_y(x) = T(x) + y$ is homotopic to $T$ by the straight line homotopy. We therefore have, by Lemma~\ref{lem:deg_homotopy}, that $\Lambda(P_{\gamma'} T_yr_{\gamma'}) = \Lambda(P_{\gamma'} Tr_{\gamma'}) = a_{\gamma'}$. By hypothesis, $a_{\gamma'} \neq 0$ for all $\Gamma' \ni \gamma' \geq \gamma'_0$ for some $\gamma'_0\in \Gamma'$.  Since $X_{\gamma'}$ is a Lefschetz space, there exists $x_{\gamma'}$ for each $\gamma'_0 \leq {\gamma'}\in \Gamma'$ such that $P_{\gamma'} T_yr_{\gamma'}(x_{\gamma'}) = x_{\gamma'}$. Since $T_y$ is approximation compact, there exists a subnet $\{x_{\gamma''}\}_{\gamma''\in \Gamma''}$ such that $\lim_{\gamma''} x_{\gamma''} = x\in C$ with $T_y(x) = x$, and therefore $T(x) + y = x$. 
	\end{proof}

	\section{Topological spaces and Lefschetz admissibility}\label{sec:metric}
	
	We now consider the case of topological spaces, and in particular metric spaces, and how to adapt the definitions and results of Section~\ref{sec:TVS} in the setting where such spaces embedd as neighborhood retracts in a TVS. This situation is closely related the one studied in \cite{Browder}. Our objective is to generalize the framework of \cite{Browder} to study fixed point properties on topological spaces that are ont necessarily metric, and that embedd in general TVSs rather than normed spaces or TVSs with completeness and/or metric assumptions. 
	
	Given a topological space $X$, we consider a triple $(X,Y,\epsilon)$ where $\epsilon : X \longrightarrow Y$ is an embedding of $X$ into a TVS $Y$. In what follows, whenever $\epsilon$ is the inclusion morphism, we will use the symbol $\iota$ instead of $\epsilon$. 
	
	\begin{definition}\label{def:approx_net_metric}
		{\rm 
				Let $(X,Y,\epsilon)$ be a triple as above, and let $\Gamma$ be a directed set. Then, we say that a family of continuous maps $\{P_\gamma : C_\gamma \longrightarrow X_\gamma\}_{\gamma\in \Gamma}$, where $C_\gamma \subset Y$ for each $\gamma$, is an {\it L-approximation net} for the triple $(X,Y,\epsilon)$, if $X_\gamma \subset Y$ is a Lefschetz space for each $\gamma\in \Gamma$, and if for every compact subset $K \subset X$ and every $U \in \mathcal U$, we can find $\gamma_0 \in \Gamma$ such that $\epsilon(K) \subset C_\gamma$ and
				\begin{eqnarray*}
						\epsilon(x) - P_\gamma(\epsilon(x)) \in U, 
				\end{eqnarray*}
				for all $x\in K$, whenever $\gamma \geq \gamma_0$. 
				If $W$ is a neighborhood of $\epsilon(X)$ and in the above each $X_\gamma \subset W$, we say that $\{P_\gamma\}_{\gamma\in \Gamma}$ is an approximation net relative to $W$. 
		}
	\end{definition}

	Similarly to Lemma~\ref{lem:approx_subnet}, we have the following. 
	\begin{lemma}\label{lem:approx_subnet_relW} 
		Let $\{P_\gamma\}_{\gamma\in \Gamma}$ denote an L-approximation net (relative to $W$) for the triple $(X,Y,\epsilon)$, and let $\{P_{\gamma'}\}_{\gamma'\in \Gamma}$ be a subnet. Then, $\{P_{\gamma'}\}_{\gamma'\in \Gamma}$ is also an L-approximation net (relative to $W$) for $(X,Y,\epsilon)$. 
	\end{lemma}

	\begin{definition}
		{\rm 
				We say that the triple $(X,Y,\epsilon)$ is Lefschetz-admissible (relative to $W$), or L-admissible (relative to $W$) for short, if it has an L-approximation net $\{P_\gamma\}_{\gamma\in \Gamma}$ (relative to $W$).  
		}
	\end{definition}

	\begin{lemma}\label{lem:restriction_L-admissible}
			Let $(X,Y,\epsilon)$ be an L-admissible triple (relative to $W$) with L-approximation net $\{P_\gamma\}_{\gamma\in \Gamma}$. If $C$ is a subset of $\epsilon(X)$, then $(C,Y,\iota)$ is an L-admissible triple (relative to $W$) with L-approximation net $\{P_\gamma\}_{\gamma\in \Gamma}$. 
	\end{lemma}
	\begin{proof}
			Let $U\in \mathcal U$, and let 
			$K$ be a compact subset of $C \subset \epsilon(X)$. Then $K = \epsilon(K')$ for some compact $K' \subset X$. Therefore, since  $(X,Y,\epsilon)$ is L-admissible, we can find $\gamma_0\in \Gamma$ such that $K = \epsilon(K') \subset C_\gamma$ and $\epsilon(x) - P_\gamma(\epsilon(x)) \in U$, for each $\gamma \geq \gamma_0$ and each $x\in K'$. Therefore $\iota(z) - P_\gamma(\iota(z)) \in U$ for each $K \ni z = \epsilon(x)$ with $x\in K'$.  
	\end{proof}

	\begin{proposition}
			Any metric space $(X,d)$ is part of an L-admissible triple $(X,Y,\epsilon)$ for some Banach space $Y$. 
	\end{proposition}
	\begin{proof}
			By the Kuratowski embedding \cite{Kuratowski,Sakai}, we can construct a Banach space $Y$ and an isometric embedding $\epsilon : X \longrightarrow Y$. The space $Y$ admits an L-approximation net (in the sense of Section~\ref{sec:TVS}) by Proposition~\ref{pro:normed_Ladmissible}. Then one verifies directly that such an L-approximation net $\{P_\gamma: Y \longrightarrow X_\gamma\}_{\gamma\in \Gamma}$ satisfies the approximation property of Definition~\ref{def:approx_net_metric} where $C_\gamma = Y$ for each $\gamma\in\Gamma$. 
	\end{proof}

	\begin{proposition}\label{pro:finite_dim_Ladmissible}
			Let $X$ be a paracompact (Hausdorff) TVS (possibly not locally convex), and let $A$ be a closed subset with finite Lebesgue covering dimension.  Then, the triple $(A, X, \iota)$ is L-admissible. 
	\end{proposition}
	\begin{proof}
				Since $X$ is paracompact and Hausdorff, it is in particular normal by Dieudonn\'e's theorem (see \cite{Dieudonne}), and therefore $A$ is a normal subspace, being a closed subset of a normal space. 
				Let $\dim A = n$ denote the Lebesgue covering dimension of $A$. For $\emptyset \neq K \subset A$ compact, and $U\in \mathcal U$ open, let $V\in \mathcal U$ be a circled neighborhood of zero such that $V + \cdots + V \subset U$, where the sum is taken $n+1$ times. Consider the covering of $K$ given by $\{x + V\}_{x\in K}$, and take a finite subcovering $\{x_i + V\}_{i=1}^r$.
				 Since $X$ is normal, and $A$ and $K$ are closed, it follows that $\dim K \leq \dim A = n$ \cite{Engelking}. Since $\{K\cap (x_i+V)\}_{i=1}^r$ is an open covering of $K$, 
				 we can therefore find a finite refinement of $\{K\cap (x_i + V)\}_{i=1}^r$ with multiplicity at most $n+1$. This means that we can find opens of $X$, $\{U_\alpha\}_{\alpha=1}^m$, such that $K\cap U_\alpha \subset K\cap (x_i+V) \subset x_i+V$ for some $i(\alpha)=1,\ldots, r$. Let $\{\eta_\alpha\}_{\alpha=1}^m\cup \{\eta_\infty\}$ be a partition of unity subordinate to the open covering $\{U_\alpha\}_{\alpha=1}^m\bigcup \{X\setminus K\}$, where each $\eta_\alpha$ has support contained in $U_\alpha$, and $\eta_\infty$ has support inside $X\setminus K$. Choose an arbitrary $y\in X$, and define the map $P_\gamma : X \longrightarrow X_\gamma$ as 
				 \begin{eqnarray*}
				 		P_\gamma(x) := \sum_{\alpha=1}^m \eta_\alpha(x)x_{i(\alpha)} + \eta_\infty(x)y, 
				 \end{eqnarray*} 
			 	where $\gamma := (K,U)$, and $X_\gamma := \overline{\conv}\{x_1,\ldots, x_r, y\}$. As shown in Section~\ref{sec:TVS}, we have that $X_\gamma$ is a Lefschetz space because it is linearly homeomorphic to a compact-convex subset of a Euclidean space, which is a Lefschetz space.  Defining the directed set $\Gamma$ of pairs $(K,U)$ with $\emptyset \neq K \subset A$ compact and $U\in \mathcal U$, ordered as in the proof of Proposition~\ref{pro:TVS_Ladmissible}, we construct a net $\{P_\gamma\}_{\gamma\in \Gamma}$. We show that this is an L-approximation net for the triple $(A,X,\iota)$. Given a nonempty compact $K \subset A$ and an open $U \in \mathcal U$, we let $\gamma_0 = (K,U)$ and consider the map $P_{\gamma_0}$. For $x\in K$, since $\dim K \leq n$, indicating by $i(\alpha)$ the index such that $K\cap U_\alpha \subset x_{i(\alpha)} + V$, and using the construction of $P_{\gamma_0}$, we find that
				\begin{equation*}
						x-P_{\gamma_0}(x) = x - \sum_{j=1}^{d} \eta_{\alpha_j}(x)x_{i(\alpha_j)} = \sum_{j=1}^{d}\eta_{\alpha_j}(x)(x-x_{i(\alpha_j)}) \in \overbrace{V + \cdots + V}^{(n+1)-{\rm times}} \subset U, 
				\end{equation*}
				where $d\leq n+1$ is the number of opens $U_\alpha$ such that $x\in U_\alpha$. 
				For each $\gamma\geq \gamma_0$, we see that $x-P_\gamma(x) \in U$ for every $x\in K$ as well by definition of ordering in $\Gamma$, therefore completing the proof. 
	\end{proof}

	In the following, we consider the case where a neighborhood retraction $r: W\longrightarrow \epsilon(X)$ is given, and the L-admissibility is relative to such $W$. Given a Lefschetz approximation net $\{P_\gamma\}_{\gamma\in \Gamma}$ for an L-admissible triple, we want to construct a subnet $\tilde \Gamma$ of $\Gamma$ such that  $X_{\tilde\gamma} \subset W$ for each $\tilde\gamma\in \tilde\Gamma$. Upon renaming $\tilde \Gamma$ by $\Gamma$, we can assume that each space $X_\gamma$ lies in $W$, and therefore we will consider this case directly. We set $r_\gamma$ to be the restriction of $r$ to $X_\gamma$ as before. 
	
	\begin{definition}\label{def:approx_comp_metric}
		{\rm 
				Let $T : \epsilon(X) \longrightarrow Y$ be a continuous map such that 
				$(T(\epsilon(X)),Y,\iota)$ is a Lefschetz admissible triple relative to $W$, and such that we have a retraction $r: W \longrightarrow \epsilon(X)$. Let $\{P_\gamma\}_{\gamma\in \Gamma}$ be a Lefschetz approximation net relative to $W$. Then, we say that $T$ is approximation compact with respect to $\{P_\gamma\}_{\gamma\in \Gamma}$ and $r$, or just approximation compact if no ambiguity arises, if the following condition is satisfied. For any subnet $\Gamma'$ of $\Gamma$, and sysetm $\{x_{\gamma'}\}_{\gamma'\in \Gamma'}$ such that $x_{\gamma'} \in X_{\gamma'}$, $Tr_{\gamma'}x_{\gamma'}\in C_{\gamma'}$, and $\lim_{\gamma'} (x_{\gamma'}-P_{\gamma'}Tr_{\gamma'}x_{\gamma'}) = 0$, there exist a subnet $\Gamma''$ of $\Gamma'$ and $x\in X$ such that $\lim_{\gamma''} x_{\gamma''} = \epsilon(x)$ and $T(\epsilon(x)) = \epsilon(x)$. 
		}
	\end{definition}

	We will also consider the maps $P^y_\gamma : C_\gamma + y \longrightarrow X_\gamma+y$ defined by shifting $P_\gamma$ as $P^y_\gamma(z) := P_\gamma(z-y) + y$, for $y\in Y$. 
	
	\begin{lemma}\label{lem:hat_approx_net}
		Let $(X,Y,\iota)$ be an L-admissible triple relative to $W$, and let $\{P_\gamma\}_{\gamma\in \Gamma}$ be an approximation net relative to $W$. Then, the maps $P^y_\gamma$ defined above constitute an approximation net relative to $W+y$, for the triple $(X+y,Y,\iota)$. 
	\end{lemma}
	\begin{proof}
		We fix $y\in Y$. 
		Define $P^y_\gamma: C^y_\gamma \longrightarrow X^y_\gamma$, where $C^y_\gamma := C_\gamma + y$, $X^y_\gamma := X_\gamma + y$, as $P^y_\gamma(w) := P_\gamma(w-y) + y$. Continuity of $P^y_\gamma$ is clear. Since $X_\gamma$ is a Lefschetz space, $X^y_\gamma$ is also Lefschetz by applying Lemma~\ref{lem:homeo_Lefschetz}. We need to show the approximation property of $\{P^y_\gamma\}$ over compact subsets of $X + y$. Let $U\in \mathcal U$, and let $K'$ be a compact subset of $X+y$. We can write $K' = K + y$ for some compact subset $K$ of $X$. Let $\gamma_0\in \Gamma$ be such that $x-P_\gamma(x)\in U$ for each $x\in K$ and each $\gamma \geq \gamma_0$. For $x'\in K'$, let $x' = x+y$ for a unique $x\in K$. Then, we have
		\begin{eqnarray*}
			x'-P^y_\gamma(x') = x+y-[P_\gamma(x+y-y)+y] = x-P_\gamma(x)\in U, 
		\end{eqnarray*}
		for each $\Gamma\ni \gamma\geq \gamma_0$. Since $X_\gamma + y \subset W + y$ for each $\gamma\in \Gamma$, it follows that $\{P^y_\gamma\}_{\gamma\in\Gamma}$ is an L-approximation net relative to $W+y$ for the triple $(X+y,Y,\iota)$. 
	\end{proof}

	The proof of the following result is similar to that of Theorem~\ref{thm:approx_compact_compact_im}, and we provide some details in this context for completeness. 
	
	\begin{theorem}\label{thm:approx_compact_compact_im_neigh_retract}
				Let $T: \epsilon(X) \longrightarrow Y$ be a compact map with $T(\epsilon(X)) \subset \epsilon(X)$, where $(X,Y,\epsilon)$ is a Lefschetz admissible triple with a neighborhood retraction $r: W\longrightarrow \epsilon(X)$, with $\epsilon(X)$ closed in $Y$. Let $\{P_\gamma\}_{\gamma\in \Gamma}$ be a Lefschetz approximation net relative to $W$. Then, $T$ is approximation compact with respect to $\{P_\gamma\}_{\gamma\in \Gamma}$ and $r$. 
	\end{theorem}
	\begin{proof}
				Applying Lemma~\ref{lem:restriction_L-admissible}, from the fact that $T(\epsilon(X)) \subset \epsilon(X)$ we have that $(T(\epsilon(X)), Y, \iota)$ is an L-admissible triple. Proceeding as in Theorem~\ref{thm:approx_compact_compact_im} for a neighborhood of zero $V$ and $U\in \mathcal U$ such that $U+U \subset V$, we can show that $\lim_{\gamma'} x_{\gamma'} - Tr_{\gamma'}x_{\gamma'} = 0$ by showing that $x_{\gamma'} - P_{\gamma'}Tr_{\gamma'}x_{\gamma'} \in U$ and $P_{\gamma'}Tr_{\gamma'}x_{\gamma'} - Tr_{\gamma'}x_{\gamma'} \in U$, where Lemma~\ref{lem:approx_subnet_relW} is used to argue that $\{P_{\gamma'}\}$ is an L-approximation net relative to $W$ for the triple $(\epsilon(X),Y,\iota)$, hence for $(T(\epsilon(X)), Y, \iota)$ again by Lemma~\ref{lem:restriction_L-admissible}. From the compactness of $\overline{T(\epsilon(X))}$ we obtain a subnet of $\{x_{\gamma'}\}$ indexed by the directed set $\Gamma''$ for which $Tr_{\gamma''}x_{\gamma''}$ converges to some $x\in \epsilon(X)$, where the fact that $\epsilon(X)$ is closed is used. Then one verifies that $\lim_{\gamma''}x_{\gamma''} = x$ and $T(x) = \lim_{\gamma''}Tr_{\gamma''}x_{\gamma''} = x$ as in the proof of Theorem~\ref{thm:approx_compact_compact_im}. 
	\end{proof}

	We want to show that in some situations it is possible to construct L-approximation nets such that $X_\gamma \subset W$ for some neighborhood rectract $\epsilon(X)$. 
	
	\begin{lemma}\label{lem:approx_W}
				Let $Y$ be a locally convex (Hausdorff) TVS, and let $\epsilon: X \longrightarrow Y$ be a closed embedding, such that $W$ is a neighborhood of $\epsilon(X)$.  Then, we can construct a directed set $\Gamma$ and an L-approximation net $\{P_\gamma\}_{\gamma\in \Gamma}$ relative to $W$.
	\end{lemma}
	\begin{proof}
				We identify $X$ with $\epsilon(X)$, and we will therefore omit $\epsilon$ throughout. Moreover, by possibly considering an open subset of $W$ that contains $\epsilon(X)$, we can assume that $W$ is open.  
				Let $U$ be an arbitrary open neighborhood of zero, and let $K$ be a nonempty compact subset of $X$. We find a convex (and circled) open $V \in \mathcal U$ such that $K + V \subset W$ and $V\subset U$ as follows. Suppose by contradiction that such open does not exist, and define the net $\mathcal U_U$ to consist of all elements of $\mathcal U$ contained in $U$, ordered by reverse inclusion, i.e. $A_1 \leq A_2 \overset{\text{def}}{\iff} A_2 \subset A_1$. Then, for each $V \in \mathcal U_U$, we can find an element $x_V \in K$ such that $x_V + V \not\subset W$. We argue that the net obtained in this procedure does not have a convergent subnet, which would be a contradiction since $K$ is compact. In fact, suppose that $\{x_V\}_{V\in \mathcal U_U}$ has a converging subnet. Let $\phi: \Gamma \longrightarrow \mathcal U_U$ be a cofinal morphism, where $\Gamma$ is a dircted set, such that $\{x_{V_\gamma}\}_{\gamma\in \Gamma}$ converges to a point $x^*\in K$, where we indicated $V_\gamma := \phi(\gamma)$. Let $U'$ be an open such that $x*+U' \subset W$. Then, we can find $U^*\in \mathcal U_U$ such that $U^*+U^* \subset U'$. By convergence of the subnet $\{x_{V_\gamma}\}_{\gamma\in \Gamma}$, we can find $V_{\gamma_0}$ with $V_{\gamma_0}\geq U^*$ such that $x_{V_{\gamma_0}}\in x^*+U^*$. Then, it follows that $x_{V_{\gamma_0}}+V_{\gamma_0} \subset x_{V_{\gamma_0}} + U^* \subset x^* + U^* + U^* \subset x^* + U' \subset W$, against the construction of the net $\{x_V\}_{V\in \mathcal U_U}$. Therefore, the net $\{x_V\}_{V\in \mathcal U_U}$ has no convergent subnet in $K$ against the fact that $K$ is compact. This shows that the claimed $V \in \mathcal U$ exists. 
				
				We now cover $K$ with the opens $x+V$, for $x\in K$, and find a finite subcover $U_i := x_i + V$ for $i=1,\ldots, n$ with $K \subset \bigcup_{i=1}^n U_i \subset W$. By construction, each $U_i$ intersects $K$ nontrivially. We now construct the nerve $L = N(\{K\cap U_i\}_{i=1}^n)$ of the covering $\{K\cap U_i\}_{i=1}^n$ of $K$, and continuously map it inside $Y$ as follows. We consider the nerve of $\{K\cap U_i\}_{i=1}^n$, so that for each $U_i$, we define a vertex $v_i$, and if $K\cap U_{i_1} \cap \cdots \cap U_{i_q} \neq \emptyset$, we construct a $(q-1)$-simplex by taking the convex linear combinations of $v_{i_j}$'s with $j=1,\ldots, q$. We define a PL map $F$ from the standard geometric realization $|L| \subset \R^n$ to $Y$ by mapping each standard Euclidean basis vector $e_i$ (which is identified with $v_i$)  to $x_i$, and affinely extending on each simplex. This procedure gives, for each compact $K$ and $V$ with $K+V \subset W$, a simplicial complex $L$, along with an image of its geometric realization into $Y$. We indicate this image by $F(\lvert L \rvert) \subset Y$. We define the directed set $\Gamma$ consisting of pairs of a nonempty compact subset of $X$, and an element of $\mathcal U$, ordered as before as $(K_1,V_1) \leq (K_2,V_2)$ if and only if $K_2\supset K_1$ and $V_2 \subset V_1$. We set $X_\gamma := F(\lvert L \rvert)$, where $\gamma = (K,U)$, and $C_\gamma := K$. We define the map $P_\gamma$ as the barycentric map of $K$ into $X_\gamma$ induced by a partition of unity subordinate to the covering $\{U_i\}_{i=1}^n$ of $K$. More precisely, if $\{\eta_i\}_{i=1}^n$ is a partition of unity on the paracompact $K$ with $\eta_i$ having support in $K\cap U_i$, we set $P_\gamma(x) := \sum_{i=1}^n \eta_i(x)x_i$. We notice that $P_\gamma(x)\in X_\gamma$ by definition of nerve, the fact that $\supp \eta_i \subset K\cap U_i$ for each $i$, and the image of the geometric realization in $Y$. 
				We need to show that $X_\gamma \subset W$, and that $X_\gamma$ is a Lefschetz space that satisfies the identity map approximation property of approximation nets with respect to $P_\gamma$.  
				
				If $x\in X_\gamma$, then $x\in F(\lvert\Delta_{i_1,\ldots, i_q}\rvert)$ for some simplex $\Delta_{i_1,\ldots, i_q}$ of $L$. Since $K\cap U_{i_1} \cap \cdots \cap U_{i_q} \neq \emptyset$, it follows that each $x_{i_j}$, $j=1,\ldots, q$, lies in $z+V$ for some $z\in K\cap U_{i_1} \cap \cdots \cap U_{i_q}$, and $x = \sum_{j=1}^q t_j x_{i_j}$, where $\sum_{j=1}^q t_j=1$. Since $V$ is convex, it follows that $x$ lies in $z+V$ as well, and since $z+V \subset W$,  it follows that $X_\gamma \subset W$. Each $X_\gamma$ is a finite polyhedron, being the image of a simplicial complex via a simplexwise affine map, by applying standard results \cite{Rou-San,Sakai}. Therefore, $X_\gamma$ is the geometric support of a finite simplicial complex and it follows that it is a Lefschetz space. Moreover, the same argument as above shows that given any nonempty compact $K \subset X$ and an open $U$, if $\gamma_0 = (K,U)$ in the definition of approximation net, then $x-P_\gamma(x) \in U$ for every $\gamma \geq \gamma_0$ and $x\in K$. 
	\end{proof}

	\begin{lemma}\label{lem:approx_W_finitedim}
			Let $Y$ be a normal (Hausdorff) TVS (not necessarily locally convex), and let $X$ be of finite Lebesgue covering dimension. Assume that $\epsilon : X \longrightarrow Y$ is a closed embedding, such that $W$ is a neighborhood of $\epsilon(X)$. Then, we can construct a directed set $\Gamma$ and an L-approximation net $\{P_\gamma\}_{\gamma\in \Gamma}$ relative to $W$.
	\end{lemma}
	\begin{proof}
			By identifying $X$ with $\epsilon(X)$, we can assume that $X\subset Y$. Since $X$ is closed, if $C \subset X$ is a closed subset of $X$, it follows that $\dim C \leq \dim X$ \cite{Engelking}. By possibly taking the interior of $W$, we can assume that $W$ is an open neighborhood of $X$. 
			We let $m := \dim X$, and we let $U$ be an arbitrary open neighborhood of zero.  Let $K\neq \emptyset$ be a compact subset of $X$, hence $\dim K \leq m$ by the previous argument. We argue that we can find $V\in \mathcal U$ such that $V + \cdots + V \subset U$, and $x+V + \cdots + V  \subset W$ for each $x\in K$, where the sum has $m+1$ terms. Define $\mathcal U_U$ to consist of the elements $V$ of $\mathcal U$ such that $V + \cdots + V \subset U$, where the sum is taken $m+1$ times, and order $\mathcal U_U$ by reverse inclusion. It follows that $\mathcal U_U$ is a directed set. As 
			in the proof of Lemma~\ref{lem:approx_W}, if the claimed $V$ did not exist, we would be able to construct a net indexed by $\mathcal U_U$ which does not admit a convergent subnet, against compactness of $K$. 
			
			Therefore, let $V \in \mathcal U$ be such that $K+V\subset W$ with the sum of $m+1$ copies of $V$ contained in $U$. We cover $K$ by elements $x+V$, with $x\in K$, and extract a finite subcover $x_i + V$ for $i=1,\ldots, r$. Since $\dim K \leq m$, we can find a refinement of $\{x_i+V\}_{i=1}^r$ with multiplicity at most $m+1$. Let $\{U_\alpha\}_{\alpha=1}^n$ be a finite collection of opens of $Y$, covering $K$, such that $\{K\cap U_\alpha\}_{\alpha=1}^n$ refines $\{x_i+V\}_{i=1}^r$ and has multiplicity at most $m+1$. In other words, $\{U_\alpha\}_{\alpha=1}^n$ has the properties that $K \subset \bigcup_{\alpha=1}^nU_\alpha$, and that for each $\alpha\in \{1,\ldots, n\}$ there exists $i(\alpha) \in \{1,\ldots,r\}$ such that $K\cap U_\alpha \subset x_{i(\alpha)} + V$, where each $x\in K$ is contained in at most $m+1$ of the $U_\alpha$'s. Upon possibly taking the intersection $U_\alpha \cap W$ for each $\alpha$, we can further assume that each open $U_\alpha$ is contained in $W$.  
			We construct the nerve of the covering and an image of its geometric realization as in Lemma~\ref{lem:approx_W} via $F$, for the covering $\{K\cap U_\alpha\}_{\alpha=1}^n$. Since $\{K\cap U_\alpha\}_{\alpha=1}^n$ has multiplicity at most $m+1$, the dimension of the simplicial complex obtained is at most $m$. Letting $\Gamma$ be defined as the directed set of pairs $(K,U)$ where $K$ is a nonempty compact set and $U \in \mathcal U$ ordered as before, it follows that we can construct $X_\gamma \subset W$ for each $\gamma = (K,U)$. The map $P_\gamma : K\longrightarrow X_\gamma$ is defined via the barycentric map associated to a partition of unity subordinate to the covering $\{K\cap U_\alpha\}_{\alpha=1}^n$ of the paracompact space $K$, as in Lemma~\ref{lem:approx_W}. 
			The space $X_\gamma$ is a Lefschetz space, as it can be shown to be a finite simplicial complex (after a suitable simplicial subdivision of a compact polyhedron). The approximation property for any nonempty compact $K\subset X$ and $U\in \mathcal U$ follows from the construction of $P_\gamma$ as in Lemma~\ref{lem:approx_W}.  
	\end{proof}

	\begin{remark}
		{\rm 
				As pointed out in Section~\ref{sec:prelim}, the notion of Lebesgue covering dimension is sometimes introduced only for normal spaces, as in \cite{Engelking}. Therefore, in that case Lemma~\ref{lem:approx_W_finitedim} has to require normality. However, as pointed out in Section~\ref{sec:prelim}, some authors also allow the usual definition to be used in arbitrary spaces (e.g. see \cite{Katetov3}). In this case, one has that the same proof of  Lemma~\ref{lem:approx_W_finitedim} is applicable without significant modifications. In the present work, we allow for the general definition, but we have emphasised the compatibility with normality in the proof of Lemma~\ref{lem:approx_W_finitedim} to have a framework that is compatible as much as possible also with the literature that considers only normal spaces, when this imposes no actual restrictions. 
		}
	\end{remark}

	The content of the previous remark is briefly formulated in the following. 
	\begin{lemma}\label{lem:no_normal}
		In the statement of Lemma~\ref{lem:approx_W_finitedim}, the assumption of being a normal space can be removed.  
	\end{lemma}

	As a consequence of Lemma~\ref{lem:no_normal}, we will use Lemma~\ref{lem:approx_W_finitedim} without exiplicitly assuming normality in what follows. As pointed out in the previous remark, those who restrict the notion of Lebesgue covering dimension to normal spaces, should interpret normality as an underlying assumption throughout. 

	\begin{definition}\label{def:gen_deg_W}
		{\rm 
			Let $T : \epsilon(X) \longrightarrow Y$ be a continuous map, where $(X,Y,\epsilon)$ is an L-admissible triple with L-approximation net $\{P_\gamma\}_{\gamma \in \Gamma}$ relative to $W$, and $r : W \longrightarrow \epsilon(X)$ is a neighborhood retraction. Then, we define the generalized Lefschetz degree, denoted by $\Lambda(T)$, as 
			\begin{eqnarray*}
				\Lambda(T) &=& \{\Lambda(P_\gamma Tr_\gamma)\ |\ \gamma\in \Gamma \}', 
			\end{eqnarray*}
			where $\Lambda(Tr_\gamma)$ is the Lefschetz degree of the map $P_\gamma Tr_\gamma : X_\gamma \longrightarrow X_\gamma$. Here the prime denotes the set of limits in $\R_\infty$ of the corresponding Lefschetz numbers along subnets of $\Gamma$ for which every composition $P_\gamma Tr_\gamma:X_\gamma\to X_\gamma$ is well defined. If no such subnet exists, we set $\Lambda(T)=\emptyset$.
		}
	\end{definition}
	
	Recall that for a map $T$, we have previously defined the $y$-translate of $T$ as $T_y(x) := T(x)+y$. The definition of generalized Lefschetz degree also applies to translated maps $T_y$. However, the following remark is relevant. 

	\begin{remark}
		{\rm 
				An important difference between the case of sets that are retractions of the whole TVS, and the neighborhood retraction considered in this section, is that the maps $P_\gamma Tr_\gamma$ might not be well defined, since we could have $Tr_\gamma(X_\gamma) \not\subset C_\gamma$. Therefore, the set $\{\Lambda(P_\gamma Tr_\gamma)\ |\ \gamma\in \Gamma \}$ might be empty, in which case the generalized degree $\Lambda(T)$ is defined to be the empty set. We will show that in some situations of interest, we are guaranteed that $\Lambda(T)$ is not empty.  
		}
	\end{remark}

	\begin{lemma}
			If $H : \epsilon(X) \times [0,1] \longrightarrow Y$ is a homotopy between two maps $F$ and $G$ as in Definition~\ref{def:gen_deg_W}, such that for each $\gamma$ we have $H(r(X_\gamma),t) \subset C_\gamma$ for every $t\in [0,1]$, then $\Lambda(F) = \Lambda(G)$. 
	\end{lemma}
	\begin{proof}
			For each $\gamma\in \Gamma$, the induced homotopy $H_\gamma : X_\gamma \times [0,1] \longrightarrow X_\gamma$ defined by $H_\gamma(x,t) = P_\gamma H(r_\gamma x,t)$ is well defined, and it is a homotopy between $P_\gamma Fr_\gamma$ and $P_\gamma Gr_\gamma$. Therefore, $\Lambda(P_\gamma Fr_\gamma) = \Lambda(P_\gamma Gr_\gamma)$ for each $\gamma\in \Gamma$. It follows that $\Lambda(F) = \Lambda(G)$. 
	\end{proof}

	\begin{theorem}\label{thm:fixed_point_NR}
		Let $T: \epsilon(X) \longrightarrow Y$ be an approximation compact map, with respect to $\{P_\gamma\}_{\gamma\in \Gamma}$ and $r: W \longrightarrow \epsilon(X)$. Then, the following facts hold. 
		\begin{itemize}
			\item[(i)] 
			If there exists a subnet $\Gamma'$ of $\Gamma$ such that $Tr_{\gamma'}(X_{\gamma'}) \subset C_{\gamma'}$ for each $\gamma' \in \Gamma'$, 
			then we have $\Lambda(T) \neq \emptyset$. 
			\item[(ii)]
			Suppose $a\in \Lambda(T)$, and let $\{a_{\gamma'}\}_{\gamma'\in \Gamma'}$ be a subnet of $\{\Lambda(P_\gamma T_yr_\gamma)\ |\ \gamma\in \Gamma \}$ such that  $\lim_{\gamma'} a_{\gamma'} = a$. Then, if $a\neq 0$ or $a = 0$ and $a_{\gamma'} \neq 0$ for all $\gamma' \geq \gamma'_0$ for some $\gamma'_0\in \Gamma'$, the fixed point equation $T(x) = x$ has a solution.  
		\end{itemize}
	\end{theorem}
	\begin{proof}
			We assume throughout that $X$ is a subset of $Y$, by identifying it with $\epsilon(X)$. We therefore will omit the embedding $\epsilon$.
		
			\noindent
			(i) Since $Tr_{\gamma'}(X_{\gamma'}) \subset C_{\gamma'}$, for each $\gamma' \in \Gamma'$ we have continuous maps $P_{\gamma'}Tr_{\gamma'}: X_{\gamma'} \longrightarrow X_{\gamma'}$, and therefore $\Lambda(P_{\gamma'}Tr_{\gamma'})$ is well defined because $X_{\gamma'}$ is a Lefschetz space. By compactness of $\R_\infty$, we have an accumulation point. 
			
			\noindent
			(ii) We have a net of maps $P_{\gamma'} Tr_{\gamma'} : X_{\gamma'} \longrightarrow X_{\gamma'}$ with nontrivial Lefschetz numbers. Then, we can find $x_{\gamma'}$ satisfying the fixed point equation $P_{\gamma'}Tr_{\gamma'}x_{\gamma'} = x_{\gamma'}$, i.e. $x_{\gamma'}- P_{\gamma'}Tr_{\gamma'}x_{\gamma'} = 0$. Since $T$ is approximation compact, it follows that there exists a subnet of $\{x_{\gamma'}\}$ converging to some $x^*\in X$, and that $T(x^*)=x^*$. 
	\end{proof}

	\begin{corollary}
			Let $T: \epsilon(X) \longrightarrow Y$ be a compact map with $\overline{T(\epsilon(X))} \subset \epsilon(X)$, where $\epsilon: X\longrightarrow Y$ is a closed embedding and we either have $Y$ locally convex, or $X$ of finite Lebesgue covering dimension. Suppose we have a retraction $r: W \longrightarrow \epsilon(X)$ where $W$ is a neighborhood of $\epsilon(X)$. Then, we have $\Lambda(T) \neq \emptyset$.
	\end{corollary}
	\begin{proof}
			As in the proof of Theorem~\ref{thm:fixed_point_NR}, by identifying $X$ with the closed subset $\epsilon(X)$ of $Y$, we can omit the symbol $\epsilon$. Let $\Gamma = \{(K,U)\ |\ \emptyset \neq K\subset X \ {\rm is\ compact,}\ U\in \mathcal U\}$ ordered as usual. We apply Lemma~\ref{lem:approx_W} or Lemma~\ref{lem:approx_W_finitedim} and construct a finite simplicial complex $X_\gamma$ along with the barycentric induced map $P_\gamma : C_\gamma \longrightarrow X_\gamma$, where $C_\gamma = K$, so that $x-P_\gamma(x) \in U$ for each $x\in K$. Therefore, $(X,Y,\iota)$ is an L-admissible triple relative to $W$ with maps $P_\gamma$. By Theorem~\ref{thm:approx_compact_compact_im_neigh_retract} we have that $T$ is approximation compact. Moreover, since $\overline{T(X)}$ is a compact subset of $X$, we have a subnet $\Gamma'$ such that $Tr_{\gamma'}(X_{\gamma'}) \subset T(X) \subset C_{\gamma'}$ for each $\gamma' \in \Gamma'$. We can therefore apply Theorem~\ref{thm:fixed_point_NR} (i). 
	\end{proof}
	
	It is in principle possible to find elements $y$ such that the conditions of Theorem~\ref{thm:fixed_point_NR}~(ii) hold for $T$ and $T_y$ simultaneously, therefore showing that the fixed point equation $x = T(x) + y$ can be solved when $x= T(x)$ can be solved. However, the main obstruction is due to the fact that the L-approximation net $\{P_\gamma\}_{\gamma\in \Gamma}$ is fixed, and it should work for both $T$ and $T_y$. We therefore consider the case where each $y$ has an associated L-approximation net $\{P^y_\gamma\}_{\gamma\in \Gamma}$, and see that the degree of $T_y$ with respect to this approximation net is equal to that of $T$ in some sufficiently regular situations. 
	
	\begin{lemma}\label{lem:retraction}
				Let $X \subset Y$ where $Y$ is a TVS. 
				Let $r: W \longrightarrow X$ be a neighborhood retraction, and let $U\in \mathcal U$ be such that $X+U \subset W$. Then there exists $U \supset V\in \mathcal U$ such that for each $y\in V$ we have a retraction $r^y: W + y \longrightarrow X$.  In particular, if $X$ is compact, such $V \in \mathcal U$ exists. 
	\end{lemma}
	\begin{proof}
			Take $V\in \mathcal U$ such that $V+V \subset U$. It is a result in uniform spaces, as discussed for example in \cite{Itzkowitz}, that if $A$ and $B$ are two uniformly separated sets in a uniform space $Y$, there exists a continuous function $\lambda: Y \longrightarrow [0,1]$ with $\lambda(A) = 0$ and $\lambda(B) = 1$. Uniformly separated here means that $E_U[A] \cap E_U[B] = \emptyset$ for some $U$, where $E_U[A]$ is the set of all points $z$ such that $(y,z) \in U$ for some $y\in A$.  With the standard uniformity of a TVS given to $Y$, we can show that for $V'\in \mathcal U$ where $V'+V' \subset V$, one has $E_{V'}[X] \cap E_{V'}[Y\setminus (X+V)] = \emptyset$, and therefore we can find a continuous function $\lambda : Y \longrightarrow [0,1]$ such that $\lambda(x) = 0$ for each $x\in X$, and  $\lambda(z) = 1$ for each $z\in Y\setminus(X+V)$. We define $r^y(z) := r(z-\lambda(z)y)$ for each $z\in W+y$, where $y\in V$. We need to verify that $r^y$ is well-defined, which means that $z-\lambda(z)y\in W$ for each $z\in W+y$. When $\lambda(z) = 1$, we have that $z-y \in W + y - y = W$. When $\lambda(z) < 1$, we have that $z \in X+V$, since $\lambda(Y\setminus(X+V)) = 1$ by construction. Then, it follows that $z-\lambda(z)y \in X +V + V\subset X+U \subset W$, where we have used the balancedness of $V\in \mathcal U$, and the fact that $V+V \subset U$ by choice of $V$. It follows that $r^y$ is well-defined. Since it is a composition of continuous maps, it follows that $r^y$ is continuous.  We just need to show that $r^y(x) = x$. This follows from the fact that $\lambda(X) = 0$ and the definition of $r^y$ by means of the retraction $r: W \longrightarrow X$. 
			
			When $X$ is compact, we can show that $U\in \mathcal U$ exists such that $X+U \subset W$ as in the proof of Lemma~\ref{lem:approx_W}, and therefore the first part of the proof guarantees the existence of $V\in \mathcal U$ such that $r^y: W+y\longrightarrow X$ is a retraction  for each $y\in V$. 
	\end{proof}

		To treat fixed point equations of type $T(x) + y = x$ with neighborhood retracts, we consider as in Section~\ref{sec:TVS} the case where we allow $T$ to be shifted in $Y$. 
	
	\begin{definition}\label{def:approx_comp_metric_Z}
		{\rm 
			Let $T : \epsilon(X) \longrightarrow Y$ be a continuous map, and let $\emptyset \neq Z \subset Y$ be such that $\epsilon(X) + y \subset W$ for each $y\in Z$. Then, we say that $T$ is approximation compact with respect to $\{P_\gamma:C^y_\gamma \longrightarrow X^y_\gamma\}_{\gamma \in \Gamma}$, $Z$ and $r : W \longrightarrow \epsilon(X)$ , or just approximation compact if no ambiguity arises, if the map $T_y(x) := T(x) + y$ is approximation compact with respect to $\{P^y_\gamma\}_{\gamma\in \Gamma}$ and $r^y$ in the sense of Definition~\ref{def:approx_comp_metric} for each $y\in Z$. 
		}
	\end{definition}

	We note that the domain and codomain of the L-approximation net $\{P_\gamma^y\}_{\gamma\in \Gamma}$, following the definitions above, are $C^y_\gamma := C_\gamma + y$ and $X^y_\gamma := X_\gamma + y$, where $P_\gamma : C_\gamma \longrightarrow X_\gamma$. Following the previous conventions, we will also set $r^y_\gamma := r^y|X^y_\gamma$. 

	\begin{lemma}\label{lem:homotopic_retractions}
				Under the same assumptions of Lemma~\ref{lem:retraction}, we have that $r$ and $r^y\circ t_y$ are homotopic, where $t_y(w) := w+y$.  
	\end{lemma}
	\begin{proof}
				Let $V\in \mathcal U$ and $\lambda: Y \longrightarrow [0,1]$ be chosen as in the proof of Lemma~\ref{lem:retraction}, and let $r^y: W+y \longrightarrow X$ be the associated retraction for a fixed $y\in V$. Then, we have $r^y\circ t_y(w) = r^y(w+y) = r(w+(1-\lambda(w+y))y)$. Define the homotopy $H: W\times [0,1] \longrightarrow X$ as $H(w,t) = r(w+t(1-\lambda(w+y))y)$. It needs to be checked that for each $t\in [0,1]$, one has $w+t(1-\lambda(w+y))y\in W$. When $w+y\notin X+V$, then $\lambda(w+y) =1$ by definition of $\lambda$, and we obtain $w+t(1-\lambda(w+y))y = w \in W$ for each $t\in [0,1]$. When $w+y\in X+V$ we have $w+t(1-\lambda(w+y))y = w+y - [1 - t(1-\lambda(w+y))]y \in X + V + V \subset X+U \subset W$, since $1 - t(1-\lambda(w+y))\in [0,1]$ for all $t\in[0,1]$, and $V$ is balanced. By construction, $H$ is continuous, and $H(w,0) = r(w)$, while $H(w,1) = r^y\circ t_y(w)$ as required.  
	\end{proof}

		From the homotopy constructed in Lemma~\ref{lem:homotopic_retractions} between $r$ and $r^y\circ t_y$, we can construct a $2$-parameter homotopy as follows. We set $\hat G(w,t,\tau) = r(w+\tau t(1-\lambda(w+ty))y)$, which can be shown to be well-defined following the same procedure used to show that $H$ in the proof of Lemma~\ref{lem:homotopic_retractions} is well-defined. With this $\hat G$, we have the diagonal homotopy $G: W\times [0,1] \longrightarrow X$ given by $G(w,t) := \hat G(w,t,t)$, which is a homotopy between $r$ and $r^y\circ t_y$. 
		
		\begin{lemma}\label{lem:homotopy}
		Let $T: \epsilon(X) \longrightarrow Y$ be an approximation compact map, with respect to $\{P_\gamma\}_{\gamma\in \Gamma}$, $\emptyset \neq Z \subset V$ and $r : W \longrightarrow \epsilon(X)$, where $V\in \mathcal U$ is such that $\epsilon(X) + V \subset W$ as in Lemma~\ref{lem:retraction}. 
		Then, for each $\gamma\in \Gamma$ such that $P_\gamma T$ is well-defined, and for each $z\in V$, we have 
		\begin{eqnarray*}
			\Lambda(P_\gamma Tr_\gamma) &=& \Lambda(P^z_\gamma T_zr^z_\gamma),  
		\end{eqnarray*}
		where $P^z_\gamma : C_\gamma + z \longrightarrow X_\gamma + z$ is defined as in Lemma~\ref{lem:hat_approx_net}, and $r^y$ is the retraction of Lemma~\ref{lem:retraction}. 
	\end{lemma}
	\begin{proof}
		Let $\gamma \in \Gamma$ and $z\in V$ be fixed. 
		We define the maps $D_t : X_\gamma \longrightarrow X_\gamma + tz$, by $D_t(x) = x+tz$. 
		We define a continuous family of maps $H_t : X\longrightarrow X_\gamma + tz$, for $t\in [0,1]$, with $H_0 = P_\gamma T$ and $H_1 = P^z_\gamma T_z$. We set $P_\gamma^{tz}(w) = P(w-tz) + tz$ for each $t\in [0,1]$, which is well-defined and continuous on $C_\gamma+tz$. Then, we define  $H_t(x) = P_\gamma^{tz} \circ T_{tz}$, and we set $\tilde H(w,t) := D^{-1}_t \circ H_t \circ G_\gamma(w,t) : X_\gamma \times [0,1] \longrightarrow X_\gamma$, where $G_\gamma := G|X_\gamma\times [0,1]$ is well-defined because $X_\gamma \subset W$. It follows that $\tilde H$ is a homotopy between $\tilde H_0 = D^{-1}_0H_0D_0=P_\gamma T r_\gamma : X_\gamma \longrightarrow X_\gamma$ and $\tilde H_1 = D^{-1}_tH_1r^z_\gamma D_1 = D^{-1}_1(P^z_\gamma T_zr^z_\gamma)D_1 : X_\gamma \longrightarrow X_\gamma$, where we have used the fact that $r_\gamma := r|X_\gamma$ and $r^z_\gamma := r^z|(X_\gamma+z)$. Therefore, we have $\Lambda(P_\gamma Tr_\gamma) = \Lambda (D^{-1}_1P^z_\gamma T_zr^z_\gamma D_1)$. However, since $D_1$ is a homeomorphism, applying the same procedure as in Lemma~\ref{lem:homeo_Lefschetz} we obtain $\Lambda(P^z_\gamma T_zr^z_\gamma) = \Lambda (D^{-1}_1P^z_\gamma T_zr^z_\gamma D_1)$, which completes the proof.   
	\end{proof}

	\begin{theorem}
			Let $T$, $\{P_\gamma\}_{\gamma\in \Gamma}$, $r: W \longrightarrow \epsilon(X)$ be as in Theorem~\ref{thm:approx_compact_compact_im_neigh_retract}, and assume that $Z \subset V \in \mathcal U$ is such that $\epsilon(X) + V \subset W$ and $T_z(\epsilon(X)) \subset \epsilon(X)$ for each $z\in Z$.  
			If $T$ satisfies (ii) in Theorem~\ref{thm:fixed_point_NR}, then $T_z$ satisfies (ii) too with respect to $\{P^z_\gamma\}_{\gamma\in \Gamma}$, and $r^z$, for each $z\in Z \cap V'$ for some $V \supset V'\in \mathcal U$. In particular, the fixed point equation $T(x) + z = x$ has a solution for each $z\in Z \cap V'$. 
	\end{theorem}
	\begin{proof}
				We show that $T_z$ is approximation compact w.r.t. $\{P^z_\gamma\}_{\gamma\in \Gamma}$ and $r^z : W + z \longrightarrow \epsilon(X)$ for each $z \in Z \cap V'$ for some $V'\in \mathcal U$ with $V' \subset V$. Lemma~\ref{lem:retraction} shows that there exists $V'\in \mathcal U$ with $V'\subset V$ such that $r^z: W+z \longrightarrow \epsilon(X)$ is a retraction for each $z\in Z \cap V'$. 
				
				Let $z\in Z\cap V'$ be fixed. 
				By Lemma~\ref{lem:hat_approx_net}, and the fact that $(T(\epsilon(X)), Y, \iota)$ is L-admissible relative to $W$ with approximation net the $\{P_\gamma\}_{\gamma\in \Gamma}$, it follows that $(T_z(\epsilon(X)), Y, \iota)$ is L-admissible relative to $W+z$ with the approximation net $\{P_\gamma^z\}_{\gamma\in \Gamma}$. Since $\overline{T(\epsilon(X))}$ is compact, it follows that $\overline{T_z(\epsilon(X))}$ is compact. 
				Then, since $T_z(\epsilon(X)) \subset \epsilon(X)$, applying the same proof as in Theorem~\ref{thm:approx_compact_compact_im_neigh_retract} it follows that $T_z$ is approximation compact with respect to $\{P^z_\gamma\}_{\gamma\in \Gamma}$ and $r^z : W + z \longrightarrow \epsilon(X)$.  
				
				Lemma~\ref{lem:homotopy} shows that $\Lambda(P_\gamma Tr_\gamma) = \Lambda(P_\gamma^zT_zr^z_\gamma)$, and therefore $T_z$ satisfies all the conditions of Theorem~\ref{thm:fixed_point_NR} (ii). In particular the fixed point equation $T_z(\epsilon(x)) = \epsilon(x)$ admits a solution. 
	\end{proof}

	\section{Some applications}\label{sec:applications}
	
	The following result can be seen as an extension of Browder's Theorem~3 and Theorem~4 in \cite{Browder} to the case of topological vector spaces without assumptions local convexity, or quasi-completeness and metrizability of the ambient TVS, and such that the embedded space is not a Lefschetz space. 
	
	For $T: X \longrightarrow X$, where $X$ embeds into $Y$ via $\epsilon$, we denote by $\hat T$ the map $\hat T: \epsilon(X) \longrightarrow Y$ induced by the embedding $\epsilon$. More specifically, we set $\hat T = \epsilon T \epsilon^{-1} : \epsilon(X) \longrightarrow \epsilon(X) \hookrightarrow Y$.  
	
	\begin{theorem}[Browder theorem in non-metrizable non-complete TVSs]
			Let $T: X \longrightarrow X$ be a compact map, with a closed embedding $\epsilon : X \longrightarrow Y$, and a neighborhood retraction $r: W \longrightarrow \epsilon(X)$. Suppose either that $Y$ is locally convex, or that $\epsilon(X)$ has finite Lebesgue covering dimension. Then, the following facts hold. 
			\begin{itemize}
				\item[(1)] 
					$\hat T$ is approximation compact w.r.t. an L-approximation net $\{P_\gamma\}_{\gamma\in \Gamma}$ and $r$;
				\item[(2)]
					$\Lambda(\hat T)\neq \emptyset$, where $\Lambda(\hat T)$ is defined via $\{P_\gamma\}_{\gamma\in \Gamma}$ and $r$ as in (1);
				\item[(3)]
					If condition (ii) in Theorem~\ref{thm:fixed_point_NR} holds for $\hat T$, where $\{P_\gamma\}_{\gamma\in \Gamma}$ and $r$ are as in (1), $T$ has a fixed point.  
			\end{itemize}
	\end{theorem}
	\begin{proof}
				Since $T$ is compact, we have that $\overline{T(X)}$ is compact in $X$, hence $\epsilon(\overline{T(X)})$ is compact in $Y$. Then, we have
				\begin{eqnarray*}
						\overline{\hat T(\epsilon(X))} = \overline{\epsilon(T\epsilon^{-1}(\epsilon(X)))} = \overline{\epsilon(T(X))} = \epsilon(\overline{T(X)}), 
				\end{eqnarray*} 
				where we used the fact that $\epsilon$ is a closed embedding in the last equality. Hence, it follows that $\hat T$ is a compact map. By Lemma~\ref{lem:approx_W} or Lemma~\ref{lem:approx_W_finitedim}, we see that we can construct an L-approximation net $\{P_\gamma\}_{\gamma\in\Gamma}$ relative to $W$, for the triple $(X,Y,\epsilon)$, which is therefore L-admissible relative to $W$. From Theorem~\ref{thm:approx_compact_compact_im_neigh_retract}, it follows that (1) holds. 
				Since $\hat T(\epsilon(X)) \subset \epsilon(X)$ has compact closure inside the closed $\epsilon(X)$, it follows that we can find a subnet $\Gamma'$ such that $P_{\gamma'}$ is well defined over $T(X)$, from which Theorem~\ref{thm:fixed_point_NR} (i) holds. 
				 In particular we have that $\Lambda(\hat T) \neq \emptyset$, which shows that (2) holds. If $\Lambda(\hat T)$ satisfies (ii) in Theorem~\ref{thm:fixed_point_NR}, then we find a fixed point $\epsilon(x) \in \epsilon (X)$ for $\hat T$. Then, from $\epsilon T \epsilon^{-1} \epsilon(x) = \epsilon(x)$ it follows that $x$ is a fixed point for $T$. 
	\end{proof}

	We now consider an application to the existence of fixed points of self-maps of compact contractible spaces that are homeomorphic to retracts of TVSs, namely Theorem~\ref{thm:compact_LC} below. It is our understanding that this result follows also from work of Askoura and Godet-Thobie \cite{Ask-GTh2} (see also the related \cite{Ask-GTh}), and Cauty \cite{Cauty}, although following a different approach. However, the stronger version of Theorem~\ref{thm:compact_LC}, namely Theorem~\ref{thm:app_compact_deg1}, seems not to follow (at least not directly) from \cite{Ask-GTh2,Cauty}. 
	
	The following result follows methods found in \cite{Sakai}. However, we have not been able to identify a reference in the literature that contains the precise statement below, which we therefore prove. 
	
	\begin{proposition}
			Let $K$ be a compact contractible subset of a Hausdorff TVS $Y$. Then, $K$ is a neighborhood retract, if and only if it is a global retract of $Y$. 
	\end{proposition}
	\begin{proof}
			The only nontrivial implication is that a neighborhood retract is a global retract. Let $r: W \longrightarrow K$ be a neighborhood retract. Upon possibly restricting to $W^\circ$, we can assume that $W$ is open. Let $a\in K$ a point such that there exists a homotopy $H : K \times [0,1] \longrightarrow K$ with $H(x,0) = a$ and $H(x,1) = x$ for each $x\in K$. We want to construct a continuous map $\lambda : Y \longrightarrow [0,1]$ with $\lambda(K) = 1$ and $\supp \lambda \subset W$. Since a Hausdorff TVS is Tikhonov, for each $x\in K$ we can find a function $u_x : Y \longrightarrow [0,1]$ such that $u_x(x) = 1$ and $u_x(Y-W) = 0$. Then, we have an open covering $\{u^{-1}_x((\frac{1}{2},1])\ |\ x\in K \}$. Let $\{u_1, \ldots, u_n\}$ be finitely many functions such that $u_i$ corresponds to $x_i\in K$ for each $i=1,\ldots, n$ and such that $\{u^{-1}_i((\frac{1}{2},1])\ |\ i=1,\ldots, n \}$ covers $K$. Let $g(x) := \max_i u_i(x)$, so that $g(x) > \frac{1}{2}$ for each $x\in K$. We define the function $\lambda$ as
			\begin{eqnarray*}
					\lambda(x) := \min \{1,\max\{0,4g(x)-1\} \}. 
			\end{eqnarray*}
			Then, by construction we have that $\lambda(x) = 1$ on $K$ and that $\supp \lambda \subset W$ as required. We define a retraction $R: Y \longrightarrow K$ by the assignment 
			\begin{eqnarray*}
					R(x) &=& \begin{cases}
							H(r(x),\lambda(x)) \ \ \ x\in W\\
							a \ \ \ x\notin W. 
					\end{cases}
			\end{eqnarray*}
			Then, $R$ is a global retraction of $Y$ onto $K$. Indeed, continuity of $R$ follows from the pasting lemma, $R(Y) \subset K$ by construction, while $R(x) = H(r(x),\lambda(1)) = H(x,1) = x$ for each $x\in K$ using the fact that $r$ is a retraction of $W$ onto $K$, the fact that $\lambda(x) = 1$ for each $x\in K$, and the definition of the homotopy $H$.
	\end{proof}

	It follows that when we consider the case of compact contractible subsets of $Y$, neighborhood retracts can be replaced by global retracts, i.e. we can use the setting of Section~\ref{sec:TVS}. Morevoer, since a (global) retract of a contractible space such as $Y$ is automatically contractible, we will consider the problem of fixed points of self-maps of compact subsets of $Y$ which are global retracts.  
	
	\begin{theorem}\label{thm:compact_LC}
			Let $T: X\longrightarrow X$ be a continuous map where $X$ is a compact subset of the locally convex TVS $Y$, with a retraction $r : Y \longrightarrow X$. Then $T$ has a fixed point. 
	\end{theorem}
	\begin{proof}
			Applying the same construction of Proposition~\ref{pro:no_paracompact}, we have that $(X,Y,\iota)$ is L-admissible with some L-approximation net $\{P_\gamma\}_{\gamma\in \Gamma}$, where each $X_\gamma$ is convex, and therefore contractible. We have that $T$ is approximation compact by Theorem~\ref{thm:approx_compact_compact_im}. Moreover, we note that by taking $K = T(X)$, we have a subnet $\Gamma'$ consisting of $\gamma$'s of type $(K',U)$, where $U\in \mathcal U$ and $K' \supset K$. Since each $X_\gamma$ is contractible, it follows that $\Lambda(P_\gamma Tr_\gamma) = 1$ for each $\gamma$ corresponding to the subnet elements, and therefore $1 \in \Lambda(T)$. Since $T$ is approximation compact, $T$ has a fixed point by Theorem~\ref{thm:fixed_point}. 
	\end{proof}
	In fact, the results of Section~\ref{sec:TVS} produce a stronger result. 
	\begin{theorem}\label{thm:app_compact_deg1}
			Let $T : C \longrightarrow X$ be approximation compact with respect to $\{P_\gamma\}_{\gamma \in \Gamma}$, $Y \subset X$, and $r : X \longrightarrow C$. Then, $\Lambda(T_y) = \{1\}$, and $T_y$ has a fixed point for every $y\in Y$. 
	\end{theorem}
	\begin{proof}
			This follows from Lemma~\ref{lem:deg_homotopy}, observing that $H_\gamma(z,t) := P_\gamma((1-t)T_y(r_\gamma(z)))$, gives a well-defined map $H_\gamma : X_\gamma \times [0,1] \longrightarrow X_\gamma$, which is a homotopy from $P_\gamma T_y r_\gamma$ and the constant map $P_\gamma(0)$. Then, since the degree of a constant map is $1$, it follows that $\Lambda(P_\gamma T_yr_\gamma) = 1$ for each $\gamma\in \Gamma$, which implies that $\Lambda(T_y) = \{1\}$. Applying Theorem~\ref{thm:fixed_point}, it follows that $T$ has a fixed point. 
	\end{proof}
		
	In a similar manner to how Browder's Theorem~3 in \cite{Browder} shows that Kinoshita's (counter)example in \cite{Kinoshita} is not a retract of $\R^3$, one can use Theorem~\ref{thm:compact_LC} to show that a compact space is not a retract of a locally convex TVS. We state this corollary as a criterion for nonexistence of retractions. 
	\begin{criterion}[No-retraction criterion]
			Let $T : X \longrightarrow X$ be a continuous map on a compact subset of a locally convex TVS $Y$. If $T$ has no fixed points, then $X$ is not a retract of $Y$. 
	\end{criterion}

\end{document}